\documentclass[12pt]{amsart}
\numberwithin{equation}{subsection}

\usepackage{amscd,amsthm,amssymb,amsfonts,amsmath,euscript}
\theoremstyle{plain}
\newtheorem{thm}[subsection]{Theorem}
\newtheorem*{thm*}{Theorem}
\newtheorem{lemma}[subsection]{Lemma}
\newtheorem{prop}[subsection]{Proposition}
\newtheorem{cor}[subsection]{Corollary}

\theoremstyle{definition}

\newtheorem{rem}[subsection]{Remark}

\theoremstyle{remark}

\newcommand{\nc}{\newcommand}

\newcommand{\ncmd}{\newcommand}
\def\makeop#1{\expandafter\def\csname#1\endcsname
  {\mathop{\rm #1}\nolimits}\ignorespaces}
\makeop{Hom}   \makeop{End}   \makeop{Aut}   \makeop{Isom}  \makeop{Pic} 
\makeop{Gal}   \makeop{ord}   \makeop{Char}  \makeop{Div}   \makeop{Lie} 
\makeop{PGL}   \makeop{Corr}  \makeop{PSL}   \makeop{sgn}   \makeop{Spf}
\makeop{Spec}  \makeop{Tr}    \makeop{Nr}    \makeop{Fr}    \makeop{disc}
\makeop{Proj}  \makeop{supp}  \makeop{ker}   \makeop{im}    \makeop{dom}
\makeop{coker} \makeop{Stab}  \makeop{SO}    \makeop{SL}    \makeop{SL}
\makeop{Cl}    \makeop{cond}  \makeop{Br}    \makeop{inv}   \makeop{rank}
\makeop{id}    \makeop{Fil}   \makeop{Frac}  \makeop{GL}    \makeop{SU}
\makeop{Trd}   \makeop{Sp}    \makeop{Tr}    \makeop{Trd}   \makeop{diag}
\makeop{Res}   \makeop{ind}   \makeop{depth} \makeop{Tr}    \makeop{st}
\makeop{Ad}    \makeop{Int}   \makeop{tr}    \makeop{Sym}   \makeop{can}
\makeop{length}\makeop{SO}    \makeop{torsion} \makeop{GSp}
\def\makebb#1{\expandafter\def
  \csname bb#1\endcsname{{\mathbb{#1}}}\ignorespaces}
\def\makebf#1{\expandafter\def\csname bf#1\endcsname{{\bf
      #1}}\ignorespaces} 
\def\makegr#1{\expandafter\def
  \csname gr#1\endcsname{{\mathfrak{#1}}}\ignorespaces}
\def\makescr#1{\expandafter\def
  \csname scr#1\endcsname{{\EuScript{#1}}}\ignorespaces}
\def\makecal#1{\expandafter\def\csname cal#1\endcsname{{\mathcal
      #1}}\ignorespaces} 
\def\doLetters#1{#1A #1B #1C #1D #1E #1F #1G #1H #1I #1J #1K #1L #1M
                 #1N #1O #1P #1Q #1R #1S #1T #1U #1V #1W #1X #1Y #1Z}
\def\doletters#1{#1a #1b #1c #1d #1e #1f #1g #1h #1i #1j #1k #1l #1m
                 #1n #1o #1p #1q #1r #1s #1t #1u #1v #1w #1x #1y #1z}
\doLetters\makebb   \doLetters\makecal  \doLetters\makebf
\doLetters\makescr 
\doletters\makebf   \doLetters\makegr   \doletters\makegr
     \def\qed{\qedmark\medbreak}%
\def\qedmark{{\enspace\vrule height 6pt width 5pt depth 1.5pt}}%
    
\def\Gm{{{\bbG}_{\rm m}}}   

\normalsize

\makeop{Bl}

\def\Spec{{\rm Spec}}

\def\Fp{{\bbF}_p}

\def\Zp{{\bbZ}_p}
\def\Qbar{\overline{\bbQ}}

\newcommand{\Z}{\mathbb Z}
\newcommand{\Q}{\mathbb Q}

\newcommand{\C}{\mathbb C}
\renewcommand{\O}{\mathcal O} % for sheaves
\newcommand{\F}{\mathbb F}

\newcommand{\<}{\langle}   %\< is not defined yet.
\renewcommand{\>}{\rangle} %\> is already defined.

\nc{\embed}{\hookrightarrow}
\newcommand{\ch}{characteristic }
\newcommand{\ac}{algebraically closed }
\newcommand{\dieu}{Dieudonn{\'e} }

\nc{\ol}{\overline}
\nc{\wt}{\widetilde}
\nc{\opp}{\mathrm{opp}}
\nc{\ul}{\underline}

\makeop{Ram}
\makeop{Rep}

\begin{document}
\title{On reduction of Hilbert-Blumenthal varieties}
\author{Chia-Fu Yu} 
\address{
Department of Mathematics \\
Columbia University \\
New York, NY 10027, USA }
\email{chiafu@math.columbia.edu}

\address{
National Center for Theoretic Sciences\\
3rd General Building, National Tsing-Hua University \\
Sect. 2 Kuang Fu Road, Hsinchu, Taiwan 30043, R.O.C.}
\email{chiafu@math.cts.nthu.edu.tw}

%\tableofcontents   % Table of Contents

\ncmd{\ur}{\mathrm{ur}}
\ncmd{\iz}{i\in \Z/f\Z}
\ncmd{\DP}{\mathrm{DP}}
\ncmd{\slope}{\mathrm{slope}}
\ncmd{\FF}{\bfF}
\ncmd{\lie}{\ul e}

%\begin{center}
%\Large\scshape
%A note on reduction of Hilbert-Blumenthal varieties
%\end{center}
%\medbreak
%\centerline{\scshape Chia-Fu Yu} 
%\centerline{\scshape  
%Columbia University and Max-Planck-Institut f{\"u}r Mathematik} 
%\centerline{\scshape 
%23 July, 2001}
%
%\footnote{partially supported by grant DMS 9800609
%from the National Science Foundation}

\begin{abstract}
Let $O_\FF$ be the ring of integers of a totally real field $\FF$ of
degree $g$. We study the reduction of the moduli
space of separably polarized abelian $O_\FF$-varieties of dimension
$g$ modulo $p$ for a fixed prime $p$.  
The invariants and related conditions for the objects in the moduli
space are discussed. We construct a scheme-theoretic stratification by
$a$-types on the Rapoport locus and study the relation with
the slope stratification. In particular, we recover the main results
of Goren and Oort [J. Alg. Geom., 2000] on the stratifications when
$p$ is unramified in $O_\FF$. We also prove the strong Grothendieck
conjecture for the moduli space in some restricted cases, particularly
when $p$ is totally ramified in $O_\FF$. 
\end{abstract}

\maketitle

\section*{Introduction}
\label{sec:00}

Let $\FF$ be a totally real number field of degree $g$ and $O_\FF$ be its
ring of integers. A Hilbert-Blumenthal variety parameterizes the 
isomorphism classes of abelian $O_\FF$-varieties of dimension $g$ with
a certain condition (and with certain additional structure). The
purpose of the condition is to exclude some bad 
points in \ch $p$ such that the integral model becomes flat. In
[R], Rapoport used the condition that the Lie algebra of the abelian
$O_\FF$-scheme over a base scheme $S$ is a locally free $O_\FF\otimes
S$-module. We will call it the Rapoport condition and the
Rapoport locus for the defined moduli space. The condition was
modified later by Deligne and Pappas [DP] in order to confirm the
properness of the compactification constructed in [R] in the case of
bad reduction. The moduli spaces defined by Deligne and Pappas are
usually referred as the Deligne-Pappas spaces. The irreducibility and
singularities of the Deligne-Pappas spaces are determined in [DP].
In the present paper we study the reduction of these moduli spaces
modulo a fixed rational prime $p$. More precisely, we consider the
moduli spaces of  abelian $O_\FF$-varieties of dimension $g$
equipped with a compatible prime-to-$p$ polarization. 
The geometry of reduction
of these moduli spaces have been studied by E. Goren and
F. Oort in [GO] when $p$ is unramified in $O_\FF$.

% When we consider the reduction of the moduli
% space modulo $p$ for a fixed prime number $p$, the condition
% defined by Deligne-Pappas [DP, (2.1.3)] can be replaced by a
% slightly weaker condition that objects are abelian $O_\FF$-schemes
% equipped with a prime-to-$p$ polarization.

Let $\scrM^{\DP}$ denote the moduli stack over $\Spec \Z_{(p)}$ of
separably polarized abelian $O_\FF$-varieties of dimension $g$.
This is a separated Deligne-Mumford algebraic stack locally of
finite type and one can identify the moduli space defined in [DP] 
as a connected component of $\scrM^{\DP}$, see [DP, 2.1]. 
We will call it the Deligne-Pappas space. Let 
$\scrM^{R}$ denote the Rapoport locus of $\scrM^{\DP}$, which
parameterizes the objects in $\scrM^{\DP}$ satisfying the Rapoport
condition. Let
$\scrM^{\DP}_p:=\scrM^{\DP}\otimes \ol{\F}_p$ and
$\scrM:=\scrM^R\otimes \ol{\F}_p$ be the reduction of $\scrM^{\DP}$
and $\scrM^R$ modulo $p$ respectively.
Let $\O:=O_\FF\otimes \Z_p$ and $k$ be an \ac field of \ch $p$.

To each abelian $O_\FF$-variety $A$ over $k$, we define two natural
invariants called the Lie type and $a$-type. Lie types are the
invariants that classify the Lie algebras $\Lie(A)$ of $A$ as
$O_\FF\otimes k$-modules, and $a$-types are those classifying the
$\alpha$-groups of $A$ [LO].
One purpose of this paper is to understand these invariants (including
the Newton polygons) using \dieu modules. Some results on the relation
among these invariants and related conditions are obtained in Section
2--3.

 A natural problem is whether these invariants arising from the
 \dieu modules in question can be realized by abelian varieties with
 the additional structure (cf.~\ref{15}). This is the integral
 analogue of a problem of Manin, see [O2]. The following theorem
 (\ref{75}) answers it affirmatively. 

\begin{thm*} {\bf 1.}
   Any quasi-polarized $p$-divisible $\O$-group $(H,\lambda,\iota)$
   (\ref{14}) over
   $k$ is isomorphic to the $p$-divisible attached to a polarized
   abelian $O_\FF$-variety.
\end{thm*}

This result implies that the strata defined by these invariants are
all non-empty. 
As an application of Thm.~1, we 
construct an explicit example of a point $s$ in the complement of the
Rapoport locus, which is both the specialization of  
a point $t_1$ in characteristic 0 and also the specialization
of an ordinary point $t_2$.  Notice that both $t_1$ and
$t_2$ are both in the Rapoport locus.
This example directly shows that the construction $\ol {\scrM^R}$ of
Rapoport is not proper over $\Spec\, \Z_{(p)}$ and a modification of the
moduli space is needed. Furthermore, it was pointed out in [DP] that the
construction of Rapoport compactifies the Deligne-Pappas space.

The proof of the algebraization theorem goes as follows. 
We first show that the formal isogeny classes are determined by the
Newton polygons, using a result of Rapoport-Zink [RZ] and of
Rapoport-Richartz [RR]. Then we prove the
weak Grothendieck conjecture. It follows that any possible Newton
polygon can be realized by an abelian variety in question. 
By the result on formal isogenies and a theorem of Tate, all the
$p$-divisible groups with additional 
structures can be realized by abelian $O_\FF$-varieties in question.

The main part of this paper is studying the strata induced from these
three invariants. The stratification by Lie types coincides with the
stratification studied by Deligne and Pappas [DP]. The largest stratum
is the Rapoport locus. We define a scheme-theoretic stratification 
by $a$-types on the Rapoport locus. 
The relation between the alpha stratification and the
slope stratification on the Rapoport locus is given. 
As the Rapoport locus is the whole space in the case of good
reduction, we recover  
the main results of E. Goren and F. Oort [GO] on the stratifications.

We now state the results. Write $\O:=O_\FF\otimes \Z_p=\oplus_{v|p} \O_v$
and let $e_v$ and $f_v$ be the ramification index and residue degree
of $v$. Let $A$ be an abelian $O_\FF$-variety over $k$ and let
$A[p^\infty]=\oplus_{v|p}H_v$ be the 
decomposition of the associated $p$-divisible group (with respect to
the $\O$-action).  
We define in (\ref{19}) the $a$-type $\ul a(H_v)$ for each component $H_v$
and put $\ul a(A):=(\ul a(H_v))_v$.
When $A$ satisfies the Rapoport condition, 
the $a$-type $\ul a(H_v)$ of each component is of the form
$(a_v^i)_{i\in \Z/f_v\Z}$, where $0\le a_v^i\le e_v$ for all $i\in
\Z/f_v\Z$. There is a natural partial order on the set of these $a$-types. 
The {\it reduced $a$-number} of $H_v$ is defined to be 
$\dim_k(\alpha_p,H_v[\pi_v])$ (\ref{19}), where $\pi_v$ is an uniformizer of
$\O_v$. 

\def\ul{\underline}
\begin{thm*}{\bf 2.}
  Let $\ul a$ be an $a$-type on $\scrM$. The closed subscheme
  $\scrM_{\ge \ul a}$ of $\scrM$ that consists of objects with $a$-type
  $\ge \ul a$ is smooth over $\Spec \ol \F_p$ 
  of pure dimension $g-|\,\ul a\,|$ (Theorem~\ref{54}). 
\end{thm*}

Let $\ul a=(\ul a_v)_v$ be an $a$-type on $\scrM$. We call 
$\ul a_v=(a_v^i)_i$ is
\emph{spaced} if $a^i_v a^{i+1}_v=0$ for all  $i\in \Z/f_v\Z$  
and $\ul a$ is \emph{spaced} if $\ul a_v$ is spaced for all $v|p$.
We put $\lambda(\ul a_v):=\max
\{|\ul b_v|; \ul b_v\le \ul a_v, \text{$\ul b_v$ is spaced} \}$
(cf. [GO, p.~112]). We refer to (\ref{110}) for the definition of the
function $s_v$, which sends certain rational numbers to possible slope
sequences at $v$.  

\begin{thm*}{\bf 3.}
  (1) 
  If $\ul a=(\ul a_v)_v$ is spaced, then the subset of $\scrM_{\ge \ul
  a}$ consisting of points whose slope sequence is $(s_v(|\ul
  a_v|))_v$ is dense in $\scrM_{\ge \ul a}$ (Theorem~\ref{68}). 

  (2)  The generic point of each irreducible component of $\scrM_{\ge \ul
  a}$ has slope sequence $\ge (s_v(\lambda(\ul a_v))_v$ (Corollary~\ref{69}).
\end{thm*}

%   . where $\lambda(\ul a):=\max
%   \{|\ul b|; \ul b\subset \ul a, \text{$\ul b$ is spaced} \}$
%   (cf. [GO, p.112]) (Corollary~\ref{69}).

\begin{thm*}{\bf 4.} 
  (1) Let $U$ be the subset of $\scrM$ consisting of points with
  reduced $a$-number one for each component $v|p$. Then 
  the strong Grothendieck conjecture holds for $U$ (Theorem~\ref{612}).

  (2) The weak Grothendieck conjecture for
  $\scrM$ holds (Corollary~\ref{615}).

  (3) The strong Grothendieck conjecture for $\scrM_p^{\DP}$  holds
  when all the residue fields of $O_\FF$ are $\F_p$ (Theorem~\ref{619}).

\end{thm*}

For the statement of the Grothendieck conjectures, we refer to
[GO. 5.1] or (\ref{113}). \\

%   Let $x:\Spec k\to \scrM$ be a geometric point of reduced
%   $a$-number one for each component $v$. Then for each  $\ul
%   m=(m_v)\in S(\ul g),(s(m_v))_v\le \slope(x)$, 
%   the reduced closed subscheme $\scrM^{\ge s(m)}_x$ of 
%   slope $\ge (s(m_v))_v$ in $\scrM_x$ is formally smooth 
%   of codimension $\sum_v \lceil m_v \rceil$ and its generic point 
%   has slope $(s(m_v))_v$ (Theorem~\ref{612}).

The methods are based on the previous works of F. Oort [O2], E. Goren
and F.~Oort [GO], and the author [Y1]. Using the explicit deformation
method developed by Norman [N], Norman-Oort [NO], T. Zink [Z2], we
construct the universal deformation of any \dieu module in the
Rapoport locus. Then we study systematically the
alpha stratification and the slope stratification on the formal
neighborhood around the point. The results above are extracted from
the formula of iterating the Frobenius map. It is
possible to extract finer information beyond the reduced $a$-number
one from our formula. 
We leave this possible generalization in the future
when the finer information becomes useful. We follow the approach of
[O2] and [Y1], which is different from that of Goren and Oort [GO]. 
Therefore, we do not repeat the computation done in loc.~cite.
F. Andreatta and E. Goren earlier obtained similar results in the case when
$p$ is totally ramified. Our work is independent from their results.

In [C], C.-L. Chai studied the combinatorial properties of Newton
points for connected reductive quasi-split algebraic groups, inspired
by earlier works of R. Kottwitz [K1], K.-Z. Li and F. Oort [LO],
M. Rapoport and M. Richartz [RR]. He gave 
a group-theoretic conjectural description of the dimensions of Newton
strata of {\it good} reduction of Shimura varieties. The main motivation
of this work is to examine whether his group-theoretic description of
Newton strata for quasi-split groups applies for the simplest case of
bad reduction. Our results support his
description even when the reductive group $G$ in question is no
longer unramified.  One may expect that Chai's description 
applies for the bad reduction case as well, under 
the assumption of the existence of good integral models of
Shimura varieties for certain special subgroups $K_p\subset G(\Q_p)$.  
We answer a question of Chai affirmatively [C, Question 7.6,
p.~984] on the dimensions of Newton strata in the Hilbert-Blumenthal
cases when $p$ is totally ramified.    

The results of this work reveal an important feature: the
stratifications on the smooth locus of bad reduction of
Hilbert-Blumenthal varieties behave very similarly as those on the good
reduction. The reader can compare the results stated above and those of
[GO] for the good reduction case. 
A fundamental question is whether this feature holds for more 
general PEL-type Shimura varieties. We have no idea but expect that
the cases of type C still hold. 

%That is the main reason that we can extend the results of
%Goren and Oort to the bad reduction case without essentially new
%techniques. 

The following is the structure of this paper.
In Section 2 we describe the structure of the \dieu module of general
polarized abelian $O_\FF$-varieties of dimension $g$. In Section 3 we classify
the formal isogeny classes explicitly.
In Section 4 we provide the normal forms of the \dieu modules in the
Rapoport locus. In Section 5 we give the natural generalization the
alpha stratification on the Rapoport locus and study its
properties. 
In Section 6 we study systematically the
alpha stratification and the slope stratification on formal
neighborhoods in the Rapoport locus, by the methods of [N], [NO], [Z2],
and [O2] as explained before. 
In Section 7 we establish a theorem of algebraization concerning
the $p$-divisible groups in question. 
In Section 8 we give the explicit example as explained before. 
In the last section we perform a
computation of the Hecke correspondence. Using this result, we
describe the singularities of the supersingular locus near the
superspecial point constructed in Section 8. \\

We should note that there is no assumption of $p$ in this paper. As
most of the time we are dealing with conditions and properties of the
associated $p$-divisible group and local properties of the moduli
spaces. It is enough to treat each component of the $p$-divisible
groups and that of the local moduli spaces by the Serre-Tate
Theorem. Without loss of generality we may assume that there is one
prime over $p$. It is clear how to state the results in this paper
without this assumption of $p$. We feel not necessary to repeat this.

We use the convenient language of algebraic stacks. The reader is free
to replace the word ``algebraic stack(s)'' by  ``scheme(s)'' 
in this paper by adding an auxiliary level structure. All schemes 
here are implicitly assumed to be locally noetherian.  \\

{\it Acknowledgments.} This paper grows out from a letter to C.-L.~Chai
responding his questions in March 2001. I am very grateful to him for his
stimulating questions and encouragements. 
I would like to thank M.~Rapoport for his comments
and inspiring questions concerning 
the Kottwitz condition. Special thanks are due to G.~Kings for useful
comments on an earlier manuscript.    
  The manuscript was completed during my stay  at the
  Max-Planck-Institut f{\"u}r Mathematik in Bonn in the summer of
  2001. I thank the Institute for the kind hospitality and the
  excellent working environment. Finally, I would like to thank the
  referee for helpful comments which improve the exposition of this
  paper significantly.  

\section{Notations, terminologies and definitions}
\label{sec:01}

\subsection{}
\label{11}
Fix a rational prime $p$. 
Let $\FF$ be a totally real field of degree $g$ and $O_\FF$ be the
ring of integers. Let $v$ be a prime of $O_\FF$ over $p$. Let $\FF_v$
denote the completion of $\FF$ at $v$ and $\O_v$ denote the ring of
integers. Denote by $e_v$ and $f_v$ the ramification index and residue
degree of $v$ respectively. Denote by $\pi_v$ a uniformizer of the
ring of integers $\O_v$ and write $g_v:=[\FF_v:\Q_p]$.

Let $v_1,\dots, v_s$ be the primes
of $O_\FF$ over $p$, $\FF_p:=\FF\otimes \Q_p=\FF_{v_1}\oplus \dots \oplus
\FF_{v_s}$, and $\O:=O_\FF\otimes \Z_p=\O_{v_1}\oplus \dots \oplus
\O_{v_s}$. Write $g_i:=g_{v_i}=[\FF_{v_i}:\Q_p]$ and $\ul g:=(g_v)_{v|p}$.

\subsection{}\label{12}
Let $k$ be a perfect field of \ch $p$. Denote by $W:=W(k)$ the ring of
Witt vectors and $B(k)$ its field of fractions. Let $\sigma$ be the
Frobenius map on $W$. 

\subsection{}\label{13}
Let $B$ be a finite dimensional semi-simple algebra over $\Q$ with
a positive involution $*$. Let $O_B$ be an order of $B$ stable under
the involution $*$. Recall that a \emph{polarized abelian
  $O_B$-variety} [Z3] is a
triple $(A,\lambda,\iota)$ where $A$ is an abelian variety,
$\iota:O_B\to \End(A)$ is a ring monomorphism and $\lambda:A\to A^t$
is a polarization satisfying the compatible condition 
$\lambda\, \iota(b^*)=\iota(b)^t \lambda$ for all $b\in O_B$.

Let $A$ be an abelian variety up to isogeny with $\iota: B\to
\End(A)$. Then the dual abelian variety $A^t$ admits a natural
$B$-action by $\iota^t(b):=\iota(b^*)^t$. The compatible condition
above is saying that the polarization $\lambda:A\to A^t$ is $O_B$-linear.

\subsection{}\label{14}
Let $B_p$ be a finite dimensional semi-simple algebra over $\Q_p$ with
an involution $*$. Let $O_p$ be an order of $B_p$ stable under
the involution $*$. A \emph{quasi-polarized $p$-divisible
$O_p$-group} is a triple $(H,\lambda,\iota)$ where $H$ is a 
$p$-divisible group, $\iota:O_p\to \End(H)$ is a ring monomorphism 
and $\lambda:H\to H^t$
is a quasi-polarization (i.e. $\lambda^t=-\lambda$) such that  
$\lambda\, \iota(b^*)=\iota(b)^t \lambda$ for all $b\in O_p$. 

For convenience, 
we also introduce the term ``quasi-polarized \dieu $O_p$-modules''
for the associated \dieu module to a quasi-polarized
$p$-divisible $O_p$-group over $k$. 
It is a \dieu module $M$ over $k$, equipped with a $W$-valued
non-degenerate alternating pairing $\<\ , \, \>$ and a $W$-linear action by
$O_p$, that satisfies the usual condition $\<ax,y\>=\<x,a^*y\>$ and
$\<Fx,y\>=\<x,Vy\>^\sigma$ for all $a\in O_p$ and $x,y\in M$.

\subsection{}\label{15}
When $B_p=B\otimes \Q_p$ and $O_p=O_B\otimes \Z_p$. We
call a quasi-polarized
$p$-divisible $O_p$-group  $(H,\lambda,\iota)$ over $k$
\emph{algebraizable} if there is a polarized abelian $O_B$-variety
$(A,\lambda_A ,\iota_A)$
over $k$ such that the associated $p$-divisible group
$(A(p),\lambda_A(p),\iota_A(p))$ is isomorphic to $(H,\lambda,\iota)$
with the additional structure over $k$.

\subsection{}
Let $S$ be a base scheme and $A$ be an abelian $O_\FF$-scheme of
relative dimension $g$ over $S$. 
Recall that the abelian $O_\FF$-scheme $A$ satisfies \emph{the Rapoport
  condition} [R] if the Lie algebra $\Lie(A)$ is a locally free
  $O_\FF\otimes_\Z \O_S$-module. Clearly this condition is local and
  open on $S$.   

We denote by $\calP(A):=\Hom_{O_\FF}(A,A^t)^{\rm sym}$ the module of
$O_\FF$-linear symmetric homomorphisms from $A$ to its dual $A^t$ over
$S$.  Notice that $\calP(A)$ is the module of global sections
of the polarization sheaf.
If $S$ is connected and $\calP(A)$ is non-zero, then $\calP(A)$ is a
rank one projective $O_\FF$-module together with a notion of
positivity. It is shown [R, Prop.~1.12] that $\calP(A)$ is non-zero
when $S$ is the spectrum of a field or an artinian ring.  

We say the abelian $O_\FF$-scheme $A$ satisfies \emph{the
  Deligne-Pappas condition} if for any connected component $S'$ of
$S$, the module $\calP(A_{S'})$ is non-zero and the induced morphism
\[ \calP(A_{S'})\otimes_{O_\FF} A_{S'}\to A^t_{S'} \]
is isomorphic.

When $S$ is a $\Z_{(p)}$-scheme, the following simpler condition plays a
similar role: the abelian $O_\FF$-scheme $A$ admits an $O_\FF$-linear 
prime-to-$p$ degree polarization. 

\subsection{} For the reader's convenience, we recall \emph{the}
Kottwitz determinant condition [K2, Sect.~5]. Let $V$ be a
one-dimensional $\FF$-vector space.    
Let $\{e_i\}$ be a $\Z$-basis of $O_\FF$ and $\ul X=(X_i)$ be some 
indeterminants. We define $f(\ul X):=\det(\sum e_i X_i;V)$. The polynomial
$f$ is in $\Z[\ul X]$, loc.~cit. 
We say that the Lie algebra $\Lie (A)$ (or the abelian $O_\FF$-scheme 
A) satisfies \emph{the Kottwitz determinant condition} if 
\[ \det(\sum e_iX_i; \Lie(A))=f(\ul X)\] 
in $\O_S[X]$. This condition does not depend on the choice of
the basis, and it is a closed condition in a family of abelian
$O_\FF$-varieties.
If $[\Lie(A)]=[O_\FF\otimes \O_S]$ in the 
Grothendieck group of $O_\FF \otimes \O_S$-modules of
finite type, then $A$ satisfies the Kottwitz determinant condition. 

% Note that this condition is
% equivalent to
%\[ \det (a; \Lie(A))=\det (a;V) \quad \text{for all $a\in O_\FF$}. \]
% This seems to weak.

\subsection{}\label{18} 
Suppose that the ground field $k$ contains the residue
fields of $O_\FF$ over $p$. We fix a place $v$ of $O_\FF$ dividing $p$. 
Let $H$ be a $p$-divisible $\O_v$-group of height $2g_v$
over $k$ and $M$
be its associated covariant \dieu module. Note that $M$ is a free
$W\otimes_{\Zp} \O_v$-module of rank two. This follows from that the
Frobenius operator $F$ induces a bijection on $M\otimes B(k)$. We
identify the set of embeddings $\Hom(\O_v^{\rm ur}, W)$ 
with $\Z/f_v \Z$ in a way that
$\sigma:i\mapsto i+1$, where $\O_v^{\rm ur}$ is the maximal {\'e}tale
extension of $\Z_p$ in $\O_v$. 
We write  $ \O_v \otimes_{\Z_p} k=\oplus_{i\in \Z/f_v\Z}\,
k[\pi_v]/(\pi_v^{e_v})$. 
 
\emph{The Lie type} of $H$ is defined to be 
\[ \lie (H):=(\{e^i_{1},e^i_{2}\})_{i\in \Z/f_v\Z} \]
if 
\[ \Lie(H)\simeq \bigoplus_{i\in Z/f_v\Z}
\left(k[\pi_v]/(\pi_v^{e^i_{1}})\oplus 
k[\pi_v]/(\pi_v^{e^i_{2}})\right) \] 
as $O_v\otimes k$-modules for some integers $e^i_1,e^i_2$.

\subsection{}\label{19}
Let $H$ and $M$ be as above. The \emph{$a$-type} of $H$ is defined to be 
\[ \ul a(H):=(\{a^i_{1},a^i_{2}\})_{i\in \Z/f_v\Z} \]
if 
\[ M/(F,V)M \simeq \bigoplus_{i\in Z/f_v\Z}
\left(k[\pi_v]/(\pi_v^{a^i_{1}})\oplus 
k[\pi_v]/(\pi_v^{a^i_{2}})\right) \] 
as $O_v\otimes k$-modules for some integers $a^i_1,a^i_2$. The usual
$a$-number is denoted by $|\ul a(H)|$, the dimensional of the $k$-vector
space $M/(F,V)M$.

If $H$ satisfies the Rapoport condition, that is, $\Lie(H)$ is a free
$\O_v\otimes_{\Zp} k$-module, then the $a$-type $\ul a(H)$ is of the
form $(\{0,a^i\})_i$ and we write $\ul a(H)=(a^i)_i$ instead. 
In this case, we define the partial order:
$(a^i)\le (b^i)$ if $a^i\le b^i$ for all $i\in \Z/f_v\Z$, and define 
$t(H):=\dim_k M/((F,V)M+\pi_v M)$, called the 
\emph{reduced (usual) $a$-number} of $H$.

\subsection{}\label{110}
The slope sequence (the Newton polygon) of $H$ we denote by  
$\slope(H)$. It is either $\{\frac{i}{g_v},\dots,
\frac{i}{g_v},\frac{g_v-i}{g_v},\dots,\frac{g_v-i}{g_v} \}$ 
for some integer $0\le i\le g_v/2$, 
or $\{\frac{1}{2},\dots,\frac{1}{2}\}$ (Lemma~\ref{32}).

We often identify a Newton polygon with its slope sequence.  
Let $S(g_v)$ denote the subset of $\Q$:
\[ S(g_v):=\{i\in \Z; 0\le i \le \frac{g_v}{2} \}\cup \{\frac{g_v}{2}\}. \]
For each $i\in S(g_v)$, we denote by $s_v(i)$ the slope sequence
$\{\frac{i}{g_v},\dots,\frac{i}{g_v},\frac{g_v-i}{g_v},
\dots,\frac{g_v-i}{g_v}\}$. The map $s_v$ identifies the set $S(g_v)$
with that of possible slope sequences arising from $p$-divisible
$\O_v$-groups. The order on $S(g_v)$ induced from $\Q$ is 
compatible with the Grothendieck specialization theorem.

There are two possible definitions for the slopes of a \dieu module. 
One uses the slopes of the $p$-divisible group (that is, the $F$-slopes of the
{\it contravariant} \dieu module, or the $V$-slopes of the covariant
one). The other just uses its $F$-slopes, no matter which \dieu
theory (covariant or contravariant) one chooses. We adopt the
latter. As \dieu modules considered in this paper are
symmetric, the choice will not effect the results.  
   
\subsection{}\label{111}
Let $A$ be an abelian $O_\FF$-variety $k$. The associated
$p$-divisible group $A(p):=A[p^\infty]$ has the decomposition 
\[ A(p)=\oplus_{v|p} H_v. \]
We define the Lie type and $a$-type of $A$ by
\[ \lie(A):=(\lie(H_v))_{v|p},\quad \ul a(A):=(\ul a(H_v))_{v|p}. \]

Set $S(\ul g):=\prod_{v|p} S(g_v)$, and for each $\ul i=(i_v)_v\in
S(\ul g)$, we write $s(\ul i):=(s_v(i_v))_v$. 
The map $s$ identifies the set $S(\ul g)$
with that of possible Newton polygons arising from abelian
$O_\FF$-varieties of dimension $g$. 
The slope sequence of $A$ is denoted by
\[ \slope(A):=(\slope(H_v))_v. \]

\subsection{}\label{112}
Let $\scrM^{\DP}$ denote the moduli stack over $\Spec \Z_{(p)}$ of
separably polarized abelian $O_\FF$-varieties of dimension $g$.
Let $\scrM^{R}$ denote the Rapoport locus of $\scrM^{\DP}$, which
parameterizes the objects in $\scrM^{\DP}$ satisfying the Rapoport
condition. 
Let $k(p)$
be the smallest finite field containing all the residue fields
$k(v_i)$. Denote by $\scrM^{\DP}_p$ the reduction
$\scrM^{\DP}\otimes_\Z k(p)$ of $\scrM^{\DP}$ and $\scrM$ the reduction
$\scrM^R \otimes_\Z k(p)$ of the Rapoport locus $\scrM^R$. 

Let $\beta$ be an admissible Newton polygon, that is $\beta\in S(\ul
g)$. We denote by $\scrM^\beta$ (resp. $\scrM^{\ge \beta}$) the reduced
algebraic substack of $\scrM$ that consists of points with Newton
polygon $\beta$ (resp. that lies over or equals $\beta$).

Let $\ul a$ be an $a$-type on $\scrM$. Let $\scrM_{\ul a}$ denote
the reduced algebraic substack of $\scrM$ that consists of points with
$a$-type $\ul a$. 
In Sect.~\ref{sec:05}, we define a closed substack, denoted by
$\scrM_{\ge a}$, of $\scrM$ so that $x\in \scrM_{\ge \ul a}({\bar k})$ 
if and only if $\ul a(x)\ge \ul a$.

\subsection{}\label{113} 
We recall the statement of the Grothendieck conjectures [GO, 5.1]
and [O2, Sect. 6]. Let $U$ be an open subset of $\scrM^{\DP}_p$. We say 
\emph{the (strong) Grothendieck conjecture holds for $U$} if for any
$x\in U$ and any admissible Newton polygon $\beta\prec \slope(x)$, 
then for any neighborhood $V$ of $x$, there is a point in $V$ whose
Newton polygon is $\beta$. 
We say \emph{the weak Grothendieck conjecture holds for
$U$} if given any chain of admissible Newton polygons
$\gamma_1\prec \gamma_2\prec\dots\prec \gamma_s$, then there exist a
chain of irreducible subschemes of $U$:  
  $V_1\supset V_2\supset \cdots\supset V_s$ such that
  $\slope(A_{\eta_i})=\gamma_i$, where $\eta_i$ is the generic
  point of $V_i$ and $A_{\eta_i}$ is the corresponding abelian variety.

\section{\dieu modules}
\label{sec:02}

\subsection{}\label{21}
Let $\FF$ be a totally real number field of degree $g$ and $O_\FF$ be its
ring of integers. To simplify notations, we will work on the case that
there is one prime of $O_\FF$ over $p$ in Section 2--6. 
The other cases can be reduced to this case if the problem is local,
as stated before. 
Let $e$ be the ramification index and $f$ be
the residue degree of this prime $v$, thus $g=f$. 
Denote by $\FF_v$ the completion of $\FF$ at $v$ and $\O$ the ring of 
integers in $\FF_v$. 
Let $\O^{\mathrm{ur}}$ denote the maximal {\'e}tale
extension of $\Z_p$ in $\O$. The ring $\O^{\mathrm{ur}}$ is isomorphic
to $W(\F_{p^f})$. Let $\pi$ be a uniformizer of $\O$. The element
$\pi$ can be chosen from $O_\FF$ and to be totally positive, by the
weak approximation. 
Let $P(T)$ be the monic irreducible polynomial of $\pi$
over $\O^{\mathrm{ur}}$.

\subsection{}\label{22}

Let $(A,\lambda,\iota)$ be a polarized abelian $O_\FF$-variety of
dimension $g$ over a
perfect field $k$ containing $\F_{p^f}$, and let $M$ be its covariant
\dieu module.    
The \dieu module $M$ is a free $\O\otimes_{\Z_p} W(k)$-module of rank two
equipped with a non-degenerate alternating pairing
\[ M\times M\to W(k) \]
such that on
which the Frobenius $F$ and Verschiebung $V$ commute with the action of
$\O$, and that $\<ax,y\>=\<x,ay\>$ and $\<Fx,y\>=\<x,Vy\>^\sigma$ 
for all $x,y\in M$ and $a\in \O$. 
We  call it briefly a \emph{quasi-polarized \dieu $\O$-module}. 

The ring $\O\otimes_{\Z_p} W(k)$ is isomorphic to  
\[ \bigoplus_{i\in \Z/f\Z} W(k)[T]/(\sigma_i(P(T)), \]
where $\sigma_i, i\in \Z/f\Z $, are embeddings of $\O^{\mathrm{ur}}$
into $W(k)$, arranged in a way that $\sigma \sigma_i=\sigma_{i+1}$. 
Set $W^i:=W(k)[T]/(\sigma_i(P(T))$ and denote again by $\pi$ the image of $T$
in $W^i$. The action of the Frobenius map $\sigma$ on
$\O\otimes_{\Z_p} W(k)$ through the right factor gives a map
$\sigma:W^i\to W^{i+1}$ which sends $a\mapsto a^\sigma$ for $a\in W(k)$ and
$\sigma(\pi)=\pi$. We also 
have $W^i\otimes_W k=k[\pi]/(\pi^e)$ and 
$\O\otimes_{\Z_p}k=\oplus_{i\in \Z/f\Z} k[\pi]/(\pi^e)$. Let 
\[ M^i:=\{x\in M\, |\, a x=\sigma_i(a) x, \ \forall\ a\in
\O^{\mathrm{ur}} \, \} \]  
be the $\sigma_i$-component of $M$, which is a free $W^i$-module of
rank two. We have the decomposition 
\[ M=M^0\oplus M^1\oplus \dots \oplus M^{f-1} \]
in which  $F:M^i\to M^{i+1}$, $V:M^{i+1}\to M^i$. 
The summands $M^i, M^j$ are orthogonal with respect to the pairing 
$\<\ ,\,\>$ for $i\not=j$. Conversely, 
a \dieu module together with such a decomposition and these properties is a
quasi-polarized \dieu $\O$-module.

\subsection{}\label{23}
 Let $\lie(A)=(\{e^i_1,e^i_2\})_i$ be the Lie type of $A$ defined in
 (\ref{18}). The invariant $\lie(A)$ has the property that
 $0\le e^i_j\le e$ for $i\in \Z/f\Z$ and $j=1,2$, $\sum_{i}e^i_1+e^i_2=g$, and
 that the invariant $e^i_1+e^i_2$ is a locally constant function in a
 family. The last one 
 follows from the fact that the $\sigma_i$-component $\Lie(A)^i$ of
 $\Lie(A)$ is a locally free sheaf. 

 In [DP], Deligne and Pappas showed that the stratum of each Lie type
 $(\{e^i_1, e^i_2\})_i$, in the Deligne-Pappas space, is a  smooth locally
 closed subscheme, and has dimension $g-2\sum_i\min\{e^i_1, e^i_2\}$
 provided the stratum is non-empty. We will see that indeed these
 strata are non-empty (\ref{75}).

\subsection{}
Let $(H,\lambda,\iota)$ be the $p$-divisible group attached to
$(A,\lambda,\iota)$. If $\lambda_1$ is another $\O$-linear
quasi-polarization on $H$, then $\lambda_1=\lambda a$ for some $a\in
E_0:=\End_\O(H)\otimes \Q_p$ with $a^*=a$, where $*$ is the involution
induced by $\lambda$. We will show that $a\in \FF_v$. Let $\lambda_0$ be
an $\O$-linear quasi-polarization of minimal degree. Then any $\O$-linear
quasi-polarization is of the form $\lambda_0 a$ for some $a\in \O$. 

The algebra $\End(H)\otimes \Q_p$ has rank $\le 4g^2$ over
$\Q_p$. Therefore, $[E_0:\FF_v]\le 4$. If $[E_0:\FF_v]=4$, then $H$ is
supersingular and $E_0$ is a quaternion algebra over $\FF_v$. In this
case, the involution is canonical. 
If $[E_0:\FF_v]=2$, then the involution on $E_0$ is non-trivial. This
follows from the non-degeneracy of the alternating pairing. 
In either case we show that the fixed elements by $*$ lie in $\FF_v$.

Similarly, we can show that given an abelian $O_\FF$-variety and let
$\lambda_0$ be an $O_\FF$-linear polarization of minimal degree at $p$,
then any $O_\FF$-linear polarization has the form $\lambda_0 a$ for some
totally positive element $a$ in $O_\FF\otimes \Z_{(p)}$. 

\def\Norm{\mathrm {Norm}}

\subsection{}\label{25}
Let $\calD^{-1}=(\pi^{{ -d}})$ be the inverse of the different of $\O$
over $\Z_p$. There is a unique $W\otimes \O$-bilinear pairing
$ (\, , ): M\times M\to W\otimes  \calD^{-1}$
such that $\<x,y\>=\Tr_{W\otimes\O /W} (x,y)$. From the uniqueness, we
have $(Fx,y)=(x,Vy)^\sigma$ for $x,y\in M$. For each $W^i$-basis
  $x^i_1,x^i_2$ of 
$M^i$, the $\pi$-adic valuation $\ord_\pi (x^i_1, x^i_2)$ is independent
of the choice of basis and the degree of the quasi-polarization is
$p^D$ [S, Chap. 1, Prop. 12], where  
\[ D=2 \sum_{i\in \Z/f\Z}\ord_p{\Norm}_{W^i/W} (\pi^{ d}
(x^i_1,x^i_2))). \]   

We can choose two $W^i$-bases $\{x^i_1,x^i_2\}, \{y^i_1, y^i_2\}$ of
$M^i$ for each $i\in \Z/f\Z$ such that
\[ Vy^{i+1}_1=\pi^{e^i_1}x^i_1, \quad Vy^{i+1}_2=\pi^{e^i_2}x^i_2. \] 
It follows from $(Vx,Vy)=p(x,y)^{\sigma^{-1}}$ that we get 
\[ \ord_\pi (y^{i+1}_1,y^{i+1}_2)=\ord_\pi
(x^i_1,x^i_2)+(e^i_1+e^i_2-e). \]  
If $\mathrm{min}_{i\in\Z/f\Z} \ord_\pi(x^i_1,x^i_2)={-d}$ (the
exponent of the inverse different), say that $i=0$
achieves the minimum, then 
\[ \ord_\pi (x^i_1,x^i_2)={-d}+\sum_{k=0}^{i-1} (e^k_1+e^k_2-e)\]
and 
\begin{equation}\label{251}
D=2 \sum_{i=1}^{f-1} \sum_{k=0}^{i-1} (e^k_1+e^k_2-e). 
\end{equation}
If $\mathrm{min}_{i\in\Z/f\Z} \ord_\pi(x^i_1,x^i_2)>{ -d}$, then we
can divide the pairing $(\, ,)$ by a power of $\pi$ such that
$\mathrm{min}_{i\in\Z/f\Z}  
\ord_\pi(x^i_1,x^i_2)={ -d}$.

\begin{lemma}\label{lm26}
  (1) There exists a number $N$, depending only on $g$, with the following
  property: for any 
  abelian $O_\FF$-variety $(A,\iota)$ of 
  dimension $g$, there is an $O_\FF$-linear polarization $\lambda$
  such that $\ord_p(\deg \lambda)\le N$.

  (2) An abelian $O_\FF$-variety $(A,\iota)$ over $k$ admits an separable
  $O_\FF$-linear polarization if and only if $\dim(A)^i$ are the same
  for $\iz$.
\end{lemma}

\begin{proof}
  (1) The statement holds as well without the assumption
  (\ref{21}), so we prove the general case instead. Let $(A,\iota)$ be
  an abelian $O_\FF$-variety. By the weak approximation,
  we can choose an $O_\FF$-linear polarization $\lambda$ such that on
  each component $H_v$ of $A(p)$ the quasi-polarization $\lambda_v$
  has the minimal degree. Then the exponent of the local degree is given in
  (\ref{251}), and let $N'$ be the sum of these local
  exponents. The number $N'$ only depends on the Lie type but not on
  the abelian variety. As there are finitely many possible Lie types
  with a fixed dimension $g$, we take $N$ to be the maximal one
  among all $N'$. 

  (2) This follows immediately from (\ref{251}).
\end{proof}

\begin{lemma}\label{lm27}
  Let $S$ be the spectrum of a artinian ring $R$ with residue field of \ch
  $p$. Let $(A,\iota)$ be an abelian $O_\FF$-scheme over $S$. Then for
  any prime $\ell\neq p$, there exists an prime-to-$\ell$ $O_\FF$-linear
  polarization on $A$.  
\end{lemma}
\begin{proof}
  We first reduce to the case that $R$ is a field $k$. Let $R$ be a
  small extension of $R_0$. Suppose there is a prime-to-$\ell$
  polarization $\lambda$ on $A\otimes_R R_0$, then $p\lambda$ extends
  over $R$. This follows from that the obstruction class lies in
  $H^2(A_k,\O_{A_k})$, which is annihilated by $p$.
  
  As the map $\Hom(A_k,B_k)\to \Hom(A_{\bar k},B_{\bar k})$ is
  co-torsion free, the map $\calP(A_k)\to \calP(A_{\bar k})$ is
  co-torsion free. It follows that $\calP(A_k)\simeq
  \calP(A_{\bar k})$. Therefore, we need to verify the case that  $k$
  is algebraically closed. 

  Let $\lambda$ be an $O_\FF$-linear polarization on $A$ and $\<\,
  ,\>$ the induced pairing on the Tate module $T_\ell(A)$. Write
  $T_\ell(A)=\oplus_{w|\ell} T_w$ into the decomposition for
  $O_\FF\otimes \Z_\ell=\oplus_{w|\ell}\O_w$. 
  Each factor $T_w$ is a free rank two
  $\O_w$-module and and the pairing $\<\, ,\>$ induces a
  non-degenerated pairing on $T_w$. Write
  $\<x,y\>=\Tr_{\O_w/\Z_\ell}(\pi_w^{-d_w}(x,y))$ for a unique lifting 
  $(\, ,):T_w\times T_w\to \O_w$ and let $c_w:=\ord_w(e_1,e_2)$, where
  $-d_w$ is the exponent of the inverse different $\calD_w^{-1}$ of
  $\O_w$ over $\Z_\ell$ and $\{e_1,e_2\}$ is a $\O_w$-basis for $T_w$.  
  By the weak approximation, we can choose a totally positive element
  $a$ in $O_\FF[\frac{1}{\ell}]$ such that $\ord_w(a)=-c_w$ for all
  $w|\ell$. Then $\lambda a$ is an $O_\FF$-linear polarization of
  degree prime-to-$\ell$. \qed
\end{proof}

\begin{prop}\label{prop28}
  Let $(A,\iota)$ be an abelian $O_\FF$-variety over $k$. Then the
  following conditions are equivalent.

  (1) $A$ satisfies the Deligne-Pappas condition. 

  (2) $A$ admits a separable $O_\FF$-linear polarization.

  (3) $[\Lie(A)]=[O_\FF\otimes k]$ in the Grothendieck group of
      $O_\FF\otimes k$-modules of finite type.   

  (4) $A$ satisfies the Kottwitz determinant condition.

  (5) $\dim_k \Lie(A)^i$ are the same for all $\iz$.
\end{prop}

\begin{proof}
  We first remark that (1)$\implies$(3) is given in [DP,
  Prop. 2.7]. The following does not depend on this result. 

  Let $\lambda\in \calP(A)$ and let $(\lambda)$ denote the submodule
  generated by $\lambda$. Then the composition $(\lambda)\otimes A\to
  \calP(A)\otimes A\to A^t$ is $\lambda$. It follows that the degree
  of the isogeny $\calP(A)\otimes A \to A^t$ divides that of $\lambda$. It
  follows from Lemma~\ref{lm27} that the isogeny $\calP(A)\otimes A \to
  A^t$ has degree a power of $p$. This shows that (1)$\iff$(2). 

  The assertion (2)$\iff$(5) is Lemma~\ref{lm26} (2). It is clear
  that (3)$\implies$(4), as the determinant function factors through
  the Grothendieck group. 

  The semi-simplification of $\Lie(A)$, as an $O_\FF\otimes k$-module, is
  $\oplus_{\iz}\, k^{d_i}$, where $d_i=\dim_k \Lie(A)^i$. It follows that 
  (3)$\iff$(5).

  If $A$ satisfies the Kottwitz determinant condition. Then $\Lie(A)$
  is a free $\O^{\rm ur}\otimes k$-module. Then (5) follows. This
  completes the proof. \qed
\end{proof}

\subsection{}\label{scheme}
Let $S$ be a $\Z_p$-scheme and let $(A,\iota)$ be an abelian
$O_\FF$-scheme over $S$. We consider the similar conditions $(1')-(5')$
for $(A,\iota)$ over $S$, where $(1'), (2')$, and $(4')$ are the same as
$(1), (2)$, and $(4)$ in (\ref{prop28}) and

\vskip 2mm

($3'$) \  Locally for the Zariski topology, $\Lie(A)$ and $O_\FF\otimes
    \O_S$ are the same in the Grothendieck group of $O_\FF\otimes
    \O_S$-modules of finite type; 
\vskip 2mm

($5'$) $\Lie(A)$ is a locally free $\O^{\rm ur}\otimes_{\Zp}
     \O_S$-module. 

\vskip 2mm
It clear that $(5')$ is an open condition and we have
    the following implications:
    $(1')\implies (2')$  and $(3')\implies (4') \implies (5')$. 

\begin{lemma}
  If $S=\Spec R$, where $R$ is a Noetherian local ring over
  $\Z_{(p)}$, then the condition $(2')$ implies the condition $(3')$.      
\end{lemma}
\begin{proof}
  If $A$ satisfies the condition $(2')$, then we have, by [DP, Prop.~2.7,
     Remark~2.8], that $2[\Lie(A)]=2[O_\FF\otimes R]$ in the
     Grothendieck group. We may assume that $R$ is complete, as $\hat
     R$ is faithfully flat over $R$. Then it suffices to check that
     $[\Lie(A)_{R_n}]=[O_\FF\otimes R_n]$ for all
     $R_n=R/\grm^{n}$. This follows immediately from the consequence
     of the Jordan-H{\"o}rder Theorem that the Grothendieck group of
     $R'$-modules (for any aritinian ring $R'$) of finite length is
     torsion-free. \qed 
\end{proof}
% If $A$ satisfies the condition (2), then we have [DP, Prop.~2.7,
%     Remark~2.8] that locally for Zariski topology, 
% $2[\Lie(A)]=2[O_\FF\otimes \O_S]$ in the Grothendieck group. It follows that
%\[ \det(\sum e_iX_i;\Lie(A))^2=\det(\sum e_i X_i;O_\FF\otimes
%    \O_S)^2,\]
%for any $\Z$-basis $\{e_i\}$ for $O_\FF$ and some indeterminants
%$X_i$. To see whether (4) is satisfied, it suffices to check that 
%\[ \det (\sum_{j=1}^{j=e}\pi^{j-1} X_j;\Lie(A)^i)=
% \det (\sum_{j=1}^{j=e}\pi^{j-1} X_j;(O_\FF\otimes
%     \O_S)^i)\]
% for all $i$, after a finite {\'e}tale base change of $S$. 
% It follows from the following 
% elementary lemma that at least for $p\neq 2$, the condition (2) implies (4).

%\begin{lemma}
% Let $R$ be a Noetherian local ring and $2$ is invertible in $R$. Let
% $f(X_1,\dots, X_n)=X_1^e+g+h$ be a polynomial in $R[\ul X]$ such
% that $g, h\in \grm_R[\ul X]$. If $f^2=(X_1^e+g)^{2}$, then $f=X_1^e+g$.
%\end{lemma}
%\begin{proof}
%   It follows from $f^2=(X_1^e+g)^{2}$ that if $h\in \grm_R^k[\ul X]$
%    then $h\in 
%   \grm_R^{k+1}[\ul X]$. Thus $h$ lies in $\cap \grm_R^n[\ul X]$, hence
%   that $h=0$. \qed
%\end{proof}

\begin{lemma}\label{lm211}
  Let $R$ be a Noetherian local ring and let $k$ be the residue
  field. Let $A$ and $B$ be abelian schemes over $R$. The
  restriction map identifies $\Hom(A,B)$ as a subgroup of
  $\Hom(A_k,B_k)$. Then for any
  prime $\ell\neq \text{\rm char}(k)$, the quotient abelian group
  $\Hom(A_k,B_k)/\Hom(A,B)$ has no $\ell$-torsions.      
\end{lemma}
\begin{proof}
  This is a slightly modification of [O4, Lemma~2.1]. We refer to loc.~cit.
  for the proof.
\end{proof}
\begin{thm}
  Let $S$ be an $\Z_{(p)}$-scheme and $(A,\iota)$ be an abelian
  $O_\FF$-scheme over $S$. If $A$ admits a separable $O_\FF$-linear
  polarization, then $A$ satisfies the Deligne-Pappas condition. 
\end{thm}
\begin{proof}
  Write $\calP'$ for the group scheme over $S$ that represents the
  functor $T\mapsto \calP(A_T)$. We first show that if $A$ satisfies
  the Deligne-Pappas condition, then $\calP'$ is a locally constant
  group scheme over $S$. We may assume that $S$ is connected and it
  suffices to show that for any connected open subset $U$ of $S$, the
  restriction map $r: \calP'(S)\to \calP'(U)$ is an isomorphism. It is
  clear that $r$ is injective. As the Deligne-Pappas condition is
  satisfied, the composition $\calP'(S)\otimes A_U\to \calP'(U)\otimes
  A_U\to A^t_U$ is isomorphic. This shows that  $\calP'(U)\otimes
  A_U\simeq A^t_U$ and $\calP'(S)\simeq \calP'(U)$.

  We now show the statement when $S=\Spec R$, where $R$ is a Noetherian
  local $\Z_{(p)}$-algebra. Let $k$ be the residue field of $R$ and we
  identify $\calP'(R)$ as a subgroup of $\calP'(k)$. It follows from
  Lemma~\ref{lm211} that if a prime $\ell\neq \text{char}(k)$, 
  then $\calP'(k)/\calP'(R)$ is $\ell$-torsion-free. 
  As $A$ admits an separable $O_\FF$-polarization, $\calP'(k)/\calP'(R)$
  is torsion free. It follows that $\calP'(R)=\calP'(k)$, hence that $A$
  satisfies the Deligne-Pappas condition.

  We now show the statement. We may first assume that $S$ is
  connected. Let $s$ be a point of $S$. Then there is a Zariski-open connected 
  neighborhood $U_s$ of $s$ such that $\calP'(U_s)\simeq 
  \calP'(\Spec O_{S,s})$, as the
  latter is generated as a $O_\FF$-module by finitely many sections. 
  It follows that $A_{U_s}$ satisfies the Deligne-Pappas condition,
  hence that $\calP'_{U_s}$ is a constant group scheme. This shows that
  $\calP'$ is constant and $\calP'(S)=\calP'(U_s)$ for any $s$. Therefore,
  $A$ satisfies the Deligne-Pappas condition.\qed   
\end{proof}

%\subsection{}
%\label{231} 
% We say an abelian $O_\FF$-variety satisfies 
%   \emph{the Deligne-Pappas condition} if it admits a separable
%    $O_\FF$-linear 
%   polarization.   
%   It is easy to see that an
%   \dieu $\O$-module can be separably quasi-polarized if and only if 
% $e^i_1+e^i_2=e$ for all $i\in \Z/f\Z$. 
% The latter is equivalent to the condition
% that the $k$-dimensions of $\Lie(A)^i$ are all the same. This is
% equivalent to the  Kottwitz determinant condition [K2, Sect. 5]. 
% It follows that the  Deligne-Pappas condition implies the Kottwitz
%   determinant condition. 
%   We will see () that any abelian $O_\FF$-variety admits an $O_\FF$-linear
%   polarization $\lambda$. The Kottwitz determinant condition then
%   implies that the polarization $\lambda$ is a power of $\pi$ times a
%   separable polarization. Therefore, these two conditions are equivalent.

\ncmd{\Def}{\mathrm{Def}}
\subsection{}\label{24}
Let $\Def[A,\iota]$ denote the equi-\ch deformation functor of the abelian
$O_\FF$-variety 
$(A,\iota)$ over $k$. It follows from crystalline theory that 
$\Def[A,\iota](k[\epsilon])=\Hom_{k\otimes O_\FF}(VM/pM,$ $M/VM)$. 
From 
$\dim_k \Hom_{k[\pi]/\pi^e} (k[\pi]/(\pi^a),k[\pi]/(\pi^b)) 
=\min\{a,b\}$, we compute that
  \begin{equation}
   \dim_k \Def [A,\iota](k[\epsilon])=
   \sum_{i\in \Z/f\Z} \sum_{1\le j,k, \le 2}\min\{e^i_j,
      e-e^i_k\}.
  \end{equation}

If $(A,\iota)$ satisfies the Deligne-Pappas condition, then 
\begin{equation}\label{eq:242}
  \dim_k \Def [A,\iota](k[\epsilon])=ef+
   2 \sum_{i\in \Z/f\Z} \min\{e^i_1, e^i_2\}. 
\end{equation}

\subsection{}
\label{214}
Assume that $(A,\lambda,\iota)$ is a \emph{separably} polarized
abelian $O_\FF$-variety, that is, it lies in the Deligne-Pappas
space. Let $\Def[A,\lambda,\iota]$ denote the equi-\ch deformation functor of
$(A,\lambda,\iota)$. We can choose a $k[\pi]/(\pi^e)$-basis $\{x^i_1,
x^i_2\}$ of $(M/pM)^i$ for each $\iz$ such that $(VM/pM)^i$ is
generated by $\{\pi^{e_1^i}x^i_1, \pi^{e_2^i}x^i_2\}$ and that 
$\<\pi^{e-1}x^i_1, x^i_2\>=1$ and $\<\pi^k x^i_1, x^i_2\>=0$ for all
$\iz$ and $0\le k<e-1$. 
The first order
universal deformation (over the deformation ring $R$) of the
abelian $O_\FF$-variety $(A,\iota)$ is given by the following data:
\[ \wt {Fil} \subset H_1^{\mathrm{cris}}(A/R),\  \wt {Fil}=\oplus 
\wt{Fil}^i, \quad  \wt{Fil}^i=<X^i_1, X^i_2>_{R[\pi]/(\pi^e)}, \]
where 
\[ X^i_1=\pi^{e^i_1}x^i_1+\sum_{j=0}^{e^i_1-1} a_j^i \pi^j
x^i_1+\sum_{j=0}^{e^i_2-1} b_j^i \pi^j x^i_2, \quad
 X^i_2=\pi^{e^i_2}x^i_2+\sum_{j=0}^{e^i_1-1} c_j^i \pi^j
x^i_1+\sum_{j=0}^{e^i_1-1} d_j^i \pi^{j+e^i_2-e^i_1} x^i_2. \]
Here we assume that $e^i_1\le e^i_2$ for simplicity. 
The condition $\wt{Fil}^i$ being isotropic is given by
 $\<X^i_1,\pi^k X_2^i\>=0$ for $0\le k\le e-1$. When $k\ge e^i_1$,
 the condition $\<X^i_1,\pi^k X_2^i\>=0$ is automatic. For $0\le k\le
 e^i_1-1$, 
 the condition $\<X^i_1,\pi^k X_2^i\>=0$ gives the equation
 $d^i_{e^i_1-k-1}+a^i_{e^i_1-k-1}=0$.
From this we conclude that
 \begin{equation}
   \label{eq:243}
    \dim_k \Def [A,\lambda,\iota](k[\epsilon])=ef+
    \sum_{i\in \Z/f\Z}  \min\{e^i_1, e^i_2\}. 
 \end{equation}
This shows that $\dim_k \Def [A,\lambda,\iota](k[\epsilon])=g$ if and
only if $(A,\lambda,\iota)$ satisfies the Rapoport condition.

% Sketch the proof for finding the basis:

% 1. We can choose a basis such that $VM/p$ is generated by $\pi^{e_1}
% x_1, \pi^{e_2} x_2$. 

% 2. From the lifting we know that $\<\pi^i x, \pi^j x\>=0$. We also
%    have $\<\pi^i x_1, \pi^j x_2\>=0$ if $i+j\ge e$. It follows  from
%    non-degeneracy that $(\<\pi^i x_1,\pi^j x_2\>)$ is
%    non-singular. Then we can find $c_j$ such that $\<\pi^i x_1, \sum
%    c_j \pi^j x_2\>=\delta_{i, e-1}$. We are done.

\begin{prop}\label{215}
The Rapoport locus is the smooth locus in the Deligne-Pappas space.
\end{prop}
\begin{proof}
  This follows from the dimension statement of [DP, Thm. 2.2] and
  (\ref{eq:243}). \qed
\end{proof}

We will compare the Kottwitz determinant condition with the
Deligne-Pappas condition in an infinitesimal neighborhood. The
following lemma will be used for Proposition~\ref{kottwitz}.

\begin{lemma}\label{216}
  (1) Let $N$ be an $n\times n$ matrix with entries in a ring $R$ such
  that the product of any two entries is zero. Let $U=U_n$, where
  $U_n$ denotes the
  following lower triangular matrix
  \[ \begin{pmatrix}
    Y_1 & 0 & \cdots & 0 \\
    Y_2 & Y_1 & \cdots & 0 \\
    \vdots & & & 0 \\
    Y_n & Y_{n-1} & \cdots & Y_1 \\  
  \end{pmatrix},\]
  for some indeterminants $Y_i$. Then
  \begin{equation}
    \label{eq:216}
    \det(U+N)=Y_1^n+\sum_{k=1}^n(\Tr_{k-1} N) U_{1,k}
  \end{equation}
  where $\Tr_k N:=n_{1,1+k}+\cdots+n_{n-k,n}$ ($\Tr_k N:=0$ if $k\ge n$) 
  and $U_{1,k}$ denotes the $(1,k)^{th}$-cofactor of $U$. 
  
  (2) If $U'=
  \begin{pmatrix}
    U_{m_1} & 0 \\ 0 & U_{m_2} \\
  \end{pmatrix}$ and $N'=
  \begin{pmatrix}
    N_{11} & N_{12} \\ N_{21} & N_{22} \\
  \end{pmatrix}$ (of same block partition), where  
  $N'$ has the same property as $N$ above and $n=m_1+m_2$, 
  then $\det(U'+N')=Y_1^n+\sum_{k=1}^n(\Tr_{k-1}
  N_{11}+\Tr_{k-1}N_{22}) U'_{1,k}$. 
\end{lemma}
\begin{proof} (1)
  Write $U=(u_1,\dots,u_n)$ and $N=(n_1,\dots, n_n)$. Then
  \begin{equation}
    \label{eq:2161}
    \begin{split}
      \det(U+N)&=(u_1+n_1)\wedge \dots \wedge (u_n+n_n)\\
               &=u_1\wedge \dots\wedge u_n+\sum_{i=1}^n u_1\wedge\dots
               u_{i-1}\wedge n_i \wedge u_{i+1}\wedge \dots\wedge u_n. 
    \end{split}
  \end{equation}
 It follows
from (\ref{eq:2161}) and the column 
expansions that $\det(U+N)=Y_1^n+\sum_{i\le j} n_{ij}U_{ij}$. From
$U_{i,j}=U_{i+1,j+1}$ ($i\le j$), we get (\ref{eq:216}).

(2) It follows from (\ref{eq:2161}) and (\ref{eq:216}) that
\begin{equation}
  \label{eq:2163}
  \begin{split}
     \det(U'+N')&=\det(U_{m_1}+N_{11})\,\det(U_{m_2}+N_{22})\\
     &=\left(Y_1^{m_1}+\sum_{k=1}^{m_1}\Tr_{k-1} N_{11}\,
     (U_{m_1})_{1,k}\right)\left(Y_1^{m_2}+\sum_{k=1}^{m_2}\Tr_{k-1} N_{22}\,
     (U_{m_2})_{1,k}\right)\\  
     &=Y_1^n+\sum_{k=1}^{n} \left( \Tr_{k-1}
      N_{11}\cdot  (U_{m_1})_{1,k}\cdot Y_1^{m_2} +\Tr_{k-1} N_{22}\cdot
     (U_{m_2})_{1,k}\cdot Y_1^{m_1} \right).  
 \end{split}
\end{equation}
Then the statement follows from $(U_{m_2})_{1,k}\cdot Y_1^{m_1}=(U_{m_1})_{1,k}\cdot Y_1^{m_2}=U'_{1,k}$.\qed
\end{proof}

\begin{prop}\label{kottwitz}
  Let $(A,\lambda,\iota)$ be a separably polarized abelian
  $O_\FF$-variety. Let $\Def[A,\iota]^K$ denote the 
  subfunctor of $\Def[A,\iota]$ that classifies the objects  satisfying the
  Kottwitz condition. Then
  $\Def[A,\iota]^K(k[\epsilon])=\Def[A,\lambda,\iota](k[\epsilon])$.
\end{prop}
\begin{proof}
  Let notations be as in (\ref{214}). 
  Let $(\wt A,\wt\iota)$ be the universal object over $R$, which is
  written as $\otimes_{\iz} R^i$. We want to compute the equations
  defined by the condition
\[ \det (\sum_{j=1}^{e}\pi^{j-1} Y_j;\Lie(\wt A)^i)=
\det (\sum_{j=1}^{e}\pi^{j-1} Y_j;(O_\FF\otimes k)^i)\]
in $\F_p[\ul Y]$, for $\iz$. As $[\wt {{Fil}^i}]=[\Lie(\wt A)^i]$ and the right hand side
is $Y_1^e$, it reduces to the condition
\[  \det (\sum_{j=1}^{e}\pi^{j-1} Y_j;
\wt {{Fil}^i})=Y_1^e\]
It suffices to compute the defining equations on each factor $R^i$.
To ease notations, we suppress the index $i$.

We have 
\[ R=k[a_i,b_j,c_k,d_\ell \ ;\,0\le i,k,\ell< e_1, 0 \le
j<e_2]/(a_i,b_j,c_k,d_\ell)^2, \]
\[ \wt {Fil}=<X_1,X_2>_{R[\pi]/(\pi^e)}=<X_1,\pi X_1,\dots,
\pi^{e_2-1}X_1, X_2, \pi X_2, \dots, \pi^{e_1-1}X_2>_R, \]
as a free $R$-module and write $\scrB$ for this $R$-basis.
We compute from (\ref{214}) that
\begin{equation}
  \label{eq:2171}
  \begin{split}
    &\text{for $e_2\le k< e$,}\quad \pi^k X_1=\sum_{j=k-e_1}^{e_2-1}
    a_{j-k+e_1} \pi^j X_1+ \sum_{j=k-e_2}^{e_1-1}
    b_{j-k+e_2} \pi^j X_2\\
    &\text{for $e_1\le k< e$,}\quad \pi^k X_2=\sum_{j=k-e_1}^{e_2-1}
    c_{j-k+e_1} \pi^j X_1+ \sum_{j=k-e_1}^{e_1-1}
    d_{j-k+e_1} \pi^j X_2\\
  \end{split}
\end{equation}
For $k\ge 1$, let $v_k$ (resp. $v_k'$) be the column vector for $\pi^{k-1}
X_1$ (resp. $\pi^{k-1} X_2$) with respect to the basis $\scrB$. The
vectors $v_k$ (resp. $v'_k$) have coordinates in $\grm_R$ except for
$k\le e_2$ (resp. for $k\le e_1$). In the exceptional case, $v_k=E_{k}$ and
$v'_k=E_{e_2+k}$, where $\{E_i\}$ is the standard basis. The
representative matrix of the endomorphism 
$\pi^{k-1}$ on the $R$-module $\wt{Fil}$ with respect to the basis
$\scrB$ 
is 
\[ [\pi^{k-1}]=(v_k,\dots, v_{k+e_2-1}, v'_k,\dots, v'_{k+e_1-1}). \]
We have 
\[ \sum_{j=1}^{e} [\pi^{j-1}]Y_j=(f_1,\dots,f_{e_2}, f'_1,
\dots,f'_{e_1}), \]
where $f_k=\sum_{j=1}^e Y_j v_{j+k-1}$ and $f'_k=\sum_{j=1}^e Y_j
v'_{j+k-1}$. Write
\[ f_k=u_k+n_k,\quad u_k=\sum_{j=1}^{e_2-k+1} Y_j v_{j+k-1}, \quad
n_k=\sum_{j=e_2-k+2}^e Y_j v_{j+k-1};\]
\[ f_k'=u'_k+n'_k, \quad u'_k=\sum_{j=1}^{e_1-k+1} Y_j v'_{j+k-1} \quad
 n'_k=\sum_{j=e_1-k+2}^e Y_j v'_{j+k-1} \]
and put
\[ U=(u_1,\dots,u_{e_2},u'_1,\dots, u'_{e_1}), \quad 
N=(n_1,\dots,n_{e_2},n'_1,\dots, n'_{e_1})=
\begin{pmatrix}
  N_{11} & N_{12} \\ N_{21} & N_{22}. 
\end{pmatrix}\] 
Then by Lemma~\ref{216} (2), we have
\begin{equation}
  \label{eq:2172}
\det (\sum_{j=1}^{e} [\pi^{j-1}]Y_j)=Y_1^e+\sum_{k=1}^{e_2}(\Tr_{k-1}
  N_{11}+\Tr_{k-1}N_{22}) U_{1,k}. 
\end{equation}
One directly computes that for $1\le k\le e_1$,
\[ \Tr_{k-1}
  N_{11}+\Tr_{k-1}N_{22}=\sum_{j=1}^{e_1-k+1}
  j(a_{e_1-k+1-j}+d_{e_1-k+1-j})Y_j. \]
As the defining equations of $\Def[A,\lambda,\iota](k[\epsilon])$ are
  $a_i+d_i=0$ for $0\le i\le e_1-1$, one has
  $\Def[A,\lambda,\iota](k[\epsilon])\subset
  \Def[A,\iota]^K(k[\epsilon])$. Conversely, let $Y_i=0$ for $i\ge 3$,
  we have 
\[ U_{1,k}=(-1)^{k+1} Y_1^{e-k} Y_2^{k-1},\]
\[ \Tr_{k-1}  N_{11}+\Tr_{k-1}N_{22}=
  (a_{e_1-k}+d_{e_1-k})Y_1+2(a_{e_1-k-1}+d_{e_1-k-1})Y_2.\]
By comparing the coefficients of (\ref{eq:2172}) with $Y^e_1$, 
we obtain the equations $a_i+d_i=0$
for $0\le i\le e_1-1$, thus $\Def[A,\lambda,\iota](k[\epsilon])\supset
  \Def[A,\iota]^K(k[\epsilon])$. This completes the proof. \qed

\end{proof}
\begin{rem}\label{26}
  (1)
  It is known [R, Prop.~1.9] that the forgetful morphism
  $\Def[A,\lambda,\iota]\to 
  \Def[A,\iota]$ is formally {\'e}tale if $(A,\lambda,\iota)$
  satisfies the Rapoport condition, and that the Rapoport locus is
  smooth. In [DP, Thm 2.2], Deligne and Pappas concluded that 
  the singular locus had codimension two. However, they actually
  showed that the complement of the Rapoport locus had codimension
  two. Proposition~\ref{215} fills the
  harmless gap of their assertion on the dimension of the singular
  locus. 
  
  (2) From (\ref{eq:242}) and (\ref{eq:243}), one can see that the
  forgetful morphism $\Def[A,\lambda,\iota]\to \Def[A,\iota]$ is not formally
  {\'e}tale anymore when $(A,\lambda,\iota)$ does not satisfy the
  Rapoport condition. In this case, the first
  order universal deformation $(\wt A, \wt \iota)$ of $(A,\iota)$ does
  not satisfy the Deligne-Pappas condition nor the Kottwitz 
  determinant condition, but the condition $(5')$ in (\ref{scheme}) still
  holds for $(\wt A, \wt \iota)$. 

  (3) It seems that the conditions $(1')-(4')$ in (\ref{scheme}) are
      equivalent when $S$ is a local Noetherian scheme. 
      Proposition~\ref{kottwitz} shows some
      evidence. One can verify this by comparing the defining
      equations from condition $(2')$ and $(4')$ in local charts, in the
      sense of Rapoport and Zink. In fact, it is not hard to verify
      the equivalence when $e=2$. However, it is quite complicated 
      in general using this method and we do not attempt to provide
      the proof here.   
\end{rem}

\subsection{}\label{27}
Let $\bfa(A)$ denote the \emph{$a$-module} of $A$, which is defined to
be the cokernel of the Frobenius map $F$ on $\Lie(A)=M/VM$:
\[ M/VM\stackrel{F}{\to} M/VM \to \bfa(A) \to 0.  \]
 Note that in the
covariant theory, the Frobenius map $F$ is induced 
from the Verschiebung morphism $V:A^{(p)}\to A$ via the
covariant functor.

Let each $\sigma_i$-component  of $\bfa(A)$ be 
\[ \bfa(A)^i\simeq k[\pi]/(\pi^{a^i_1}) \oplus k[\pi]/(\pi^{a^i_2}) \] 
for some integers $\{a^i_1,a^i_2\}$.
We define the \emph{$a$-type $\ul a(A)$ of $A$} to be the invariant 
$\ul a(A)=(\{a^i_1,a^i_2\})_i$, (\ref{19}). Let $|\ul a(A)|$ denote
the total $a$-number 
of $A$, which is the dimension of the $k$-vector space
$\Hom_k(\alpha_p, A^t)$. If the abelian $O_\FF$-variety $(A,\iota)$
satisfies the Rapoport condition, then $\ul a(A)$ is of the form
$(\{0,a^i\})_i$ and we write $\ul a(A)=(a^i)_i$ instead.

\begin{lemma}
  The $a$-module $\bfa(A)$ is canonically isomorphic to the $k$-linear
  dual $\Hom_k(\alpha_p, A^t)^*$ of $\Hom_k(\alpha_p, A^t)$ as
  $k\otimes_\Z O_\FF$-modules. 
\end{lemma}

\begin{proof}
  Let $M^*$ denote the \emph{contravariant} \dieu functor. Then we
  have 
  \begin{equation*}
  \begin{split} 
  \Hom_k(\alpha_p, A^t)
  & =\Hom_{W[F,V]}(M^*(A^t),k)\simeq\Hom_{W[F,V]}(M(A),k) \\
  & =\Hom_k(M(A)/(F,V)M(A),k)=\bfa(A)^*. \text{\qed} 
  \end{split}
  \end{equation*} 
\end{proof}

%\subsection{}
%\label{2701}
% Note that $A$ is ordinary if and only if $|a(A)|=0$. It is easy to see
% that the latter condition implies that $\{e^i_1, e^i_2\}=\{0,e\}$ for
% all $i\in \Z/f\Z$. Therefore, ordinary points satisfy 
% the Rapoport condition automatically. As non-ordinary points do not
% specialize to the boundary of Deligne-Pappas space, the complement of
% the Rapoport locus is in the interior of the space. 
% Therefore, the compactification constructed by Rapoport provides 
% that of the Deligne-Pappas space [DP, Introduction].   

\subsection{}
\label{271}
Let $\ul a(A)=(\{a^i_1,a^i_2\})_i$, with $a^i_1\le a^i_2$, be the $a$-type of
$A$ and $\lie(A)=(\{e^i_1,e^i_2\})_i$ be the Lie type. 
It follows from the elementary divisor lemma that there are two $W^i$-bases
$\{x^i_1,x^i_2\}, \{y^i_1,y^i_2\}$ of $M^i$ such that 
\[ VM^{i+1}=< \pi^{ e^{i}_1} x^i_1, \pi^{ e^{i}_2} x^i_2>
  \text{ and }\
FM^{i-1}=<\pi^{e- {e^{i-1}_1}} y^i_1, \pi^{e- {e^{i-1}_2}} y^i_2>.\] 
We may assume that $e^i_1\le e^i_2$ hence $e-e^{i-1}_2\le e-e^{i-1}_1$
for each $i$.

If $e^i_1\le e- {e^{i-1}_2}$, write
$y^i_1=ax^i_1+bx^i_2, y^i_2=cx^i_1+dx^i_2$, then 
\[FM^{i-1}+VM^{i+1}=<\pi^{e^i_1}x^i_1, \pi^{e^i_2}x^i_2, \pi^{e-
  e^{i-1}_1}bx^i_2, \pi^{\ e-e^{i-1}_2}dx^i_2>. \]  
Note that one of $b$ and $d$ is a unit. From this we obtain a bound 
for $a$-types: $a^i_1=e^i_1$ and $\min\{e^i_2,
e-{e^{i-1}_2} \}\le a^i_2\le \min\{e^i_2,e-{e^{i-1}_1}\} $ when
$e^i_1\le e- {e^{i-1}_2}$. Conversely
if $e-{e^{i-1}_2}\le e^i_1$
then we have $a^i_1=e-{e^{i-1}_2}$ and $\min\{
e-{e^{i-1}_1},e^i_1 \}\le a^i_2\le \min\{e-{e^{i-1}_1, e^i_2}\}$.

\subsection{}\label{222}
If $A$ is superspecial, that is $|\ul a(A)|=g=ef$. We know that $FM=VM$, in
other words that $FM^{i-1}=VM^{i+1}$ for all $i\in \Z/f\Z$. This gives
$\{e^i_1,e^i_2\}=\{e-{e^{i-1}_1},e-{e^{i-1}_2}\}$ and
$\{a^i_1,a^i_2\}=\{e^i_1,e^i_2\}$ for all $i$. There
are two possibilities:

1. If $f$ is odd, then $\ul a(A)=\lie(A)=(\{e_1,e_2\})_i$ for some nonnegative
   integers $e_1$ and $e_2$ with $e_1+e_2=e$.
 
2. If $f$ is even, then there are two nonnegative integers $e_1$, $e_2$, with $0\le
   e_1,e_2\le e$, so that
   $\ul a(A)=\lie(A)=(\{e_1,e_2\},\{e-e_1,e-e_2\},\{e_1,e_2\},\dots,
   \{e-e_1,e-e_2\})$.  
\begin{lemma}
  Let $\ul a(A)=(\{a^i_1, a^i_2\})_i$ and $\lie(A)=(\{e^i_1,
  e^i_2\})_i$. Then $\lie(A^t)=(\{e-e^i_1,e-e^i_2\})_i$ and the
  $a$-type $(\{b^i_1, b^i_2\})_i$ of the dual abelian variety $A^t$
  is given as follows: $b^i_1=\min\{e-e^i_1,e-e^i_2,{e^{i-1}_1},
  {e^{i-1}_2}\}$ and $b^i_1+b^i_2=(a^i_1+a^i_2)+(e^{i-1}_1+
  e^{i-1}_2)-(e^i_1+e^i_2)$.
\end{lemma}
\begin{proof}
  Write $\ol M=M/pM$ and $\ol {M^t}=M^t/pM^t$, where $M^t$ is the \dieu
  module of $A^t$. From the definition of $\lie(A)=(\{e^i_1,
  e^i_2\})_i$, one concludes that $V\ol{ M^{i+1}}$ is isomorphic to
  $k[\pi]/(\pi^{e-e^i_1})\oplus k[\pi]/(\pi^{e-e^i_2})$. The perfect
  pairing $\ol M\times \ol {M^t}\to k$ induces a perfect pairing 
\[ V\ol {M^{i+1}}\times \ol {M^{t,i}}/ V\ol {M^{t, i+1}} \to k, \] 
thus $\lie(A^t)=(\{e-e^i_1,e-e^i_2\})_i$.

Let $T^i:=F\ol {M^{i-1}}\cap V \ol {M^{i+1}}$ in $\ol {M^i}$. We have
$\dim_k T^i=\dim_k F\ol {M^{i-1}}+\dim _k  V \ol {M^{i+1}}-\dim_k (F\ol
{M^{i-1}}+V \ol {M^{i+1}})=(a^i_1+a^i_2)+(e^{i-1}_1+
  e^{i-1}_2)-(e^i_1+e^i_2)$. The perfect pairing $\ol M\times \ol
  {M^t}\to k$ induces a perfect pairing  
  \[ T^i\times \left(M^t/ (F,V)M^t\right)^i\to k, \]
  thus $b^i_1+b^i_2=(a^i_1+a^i_2)+(e^{i-1}_1+
  e^{i-1}_2)-(e^i_1+e^i_2)$. The assertion 
  $b^i_1=\min\{e-e^i_1,e-e^i_2,{e^{i-1}_1}, 
  {e^{i-1}_2}\}$ is obtained from (\ref{271}) and
  $\lie(A^t)=(\{e-e^i_1,e-e^i_2\})_i$.\qed 
\end{proof}

\section{Formal isogeny classes}
\label{sec:03}

In the rest of this paper, $k$ denotes an \ac field of \ch $p>0$.

\begin{lemma}\label{32}
  Let $M$ be a quasi-polarized \dieu $\O$-module. Then the slope sequence
  $\mathrm{slope}(M)$ of $M$ is either $\{\frac{i}{g},\dots,
\frac{i}{g},\frac{g-i}{g},\dots,\frac{g-i}{g} \}$ for some integer $0\le i\le
  g/2$, or $\{\frac{1}{2},\dots,\frac{1}{2}\}$.  
\end{lemma}
\begin{proof}
  Suppose that $M$ is not supersingular. The F-isocrystal $M\otimes
   B(k)$ contains 
  $M_{a,b}^r\oplus M_{b,a}^r$, where $M_{a,b}$ is the simple
  F-isocrystal of single slope $a/(a+b)$ for some integers $a,b$. We want to show that
  $2(a+b)r=2g$. We first have an embedding $\FF_v\to
  \End (M_{a,b}^r)=M_r(D_{a/(a+b)})$ and any maximal commutative
  subalgebra of the latter 
  has $\Q_p$-dimension $r(a+b)$, so $g |\, r(a+b)$. On the other hand, we
  have $2g \ge 2r(a+b)$, from $M\otimes B(k)\supset M_{a,b}^r\oplus
  M_{b,a}^r$. Hence $M\otimes B(k)=M_{a,b}^r\oplus M_{b,a}^r$. This
  completes the proof.  \qed
\end{proof}

\subsection{}
\label{321}
Let $S(g)$ denote the subset of $\Q$ which parameterizes possible
slope sequences arising from abelian $O_\FF$-varieties
\[ S(g):=\{i\in \Z; 0\le i \le \frac{g}{2} \}\cup \{\frac{g}{2}\}. \]
For each $i\in S(g)$, we denote by $s(i)$ the slope sequence
$\{\frac{i}{g},\dots, \frac{i}{g},\frac{g-i}{g},\dots,\frac{g-i}{g}
\}$. The (linear) order on $S(g)$ induced from $\Q$ is compatible with
the Grothendieck specialization theorem.
   
\subsection{Example.} \label{33}
Let $M=M^0\oplus \dots \oplus M^{f-1}$ be a \dieu $\O$-module, where
$M^i$ is a free $W^i$-module of rank two generated by $x^i_1, x^i_2$
with $x^i_2\in VM^{i+1}$. Let 
\begin{equation}
  \begin{split}
    F x^i_1&=a^{i+1}_{11}x^{i+1}_1+a^{i+1}_{12}x^{i+1}_2 \\
    F x^i_2&=pa^{i+1}_{21}x^{i+1}_1+pa^{i+1}_{22}x^{i+1}_2 \\
  \end{split}
\end{equation}
be the action of the Frobenius $F$. Then 
\begin{equation}
  \begin{split}
    F^f x^0_1&=\alpha x^{0}_1+\beta x^{0}_2 \\
    F^f x^0_2&=\gamma x^{0}_1+\delta x^{0}_2, \\
  \end{split}
\end{equation}
where
\[ 
\begin{pmatrix}
  \alpha & \beta \\ \gamma & \delta \\
\end{pmatrix}=A_1^{(f-1)}\cdots A^{(1)}_{f-1} A_f, \quad
A_i:= 
\begin{pmatrix}
   a^{i}_{11} &  a^{i}_{12} \\
   pa^{i}_{21} & pa^{i}_{22} \\
\end{pmatrix}, \ A_f=A_0. 
\]
Here we write $A^{(n)}$ for $A^{\sigma^n}$. Note that $M$ satisfies
the Rapoport condition.

Let $a, b\in \Z, a+b=g,  0 \le a \le b$. Write $a=de+r, 0 \le
r<e$. Take 
\[ A_i=
\begin{pmatrix}
  0 & 1 \\ -p & 0 \\
\end{pmatrix} \ \text{for $1\le i \le 2d$}, \quad
A_i=
\begin{pmatrix}
  1 & 1 \\ -p & 0 
\end{pmatrix} \  \text{for $2d< i<f$}, \ \text{and}\ 
A_f=
\begin{pmatrix}
  \pi^r & 1 \\ -p & 0 \\
\end{pmatrix}. \]
Then we have 
\[ 
\begin{pmatrix}
  \alpha & \beta \\ \gamma & \delta \\
\end{pmatrix}=(-1)^d p^d
\begin{pmatrix}
  \alpha' & \beta' \\ \gamma' & \delta' \\
\end{pmatrix}, \quad 
\begin{pmatrix}
  \alpha' & \beta' \\ \gamma' & \delta' \\
\end{pmatrix}=\prod_{i-2d+1}^f A_i\equiv 
\begin{pmatrix}
  \pi^r & 1 \\ 0 & 0 \\
\end{pmatrix}
\mod p 
\]
Note that $A_i^\sigma=A_i$ for all $i\in \Z/f\Z$. The \ch
polynomial of $F^f$ is $X^2-(\alpha+\delta)X+p^f$ and
$\ord_p(\alpha+\delta)=\frac{a}{e}$. Therefore
$\mathrm{slope}(M)=\{\frac{a}{g},\dots, \frac{a}{g},\frac{b}{g},\dots,
\frac{b}{g}\}$. 

\begin{rem}
   The example constructed above is non-trivial, due to the following
   constraint. 
   Let $M=M_{a,b}\oplus M_{b,a},\, a+b=g,\, a\neq b$ be a \dieu
   $\O$-module which satisfies the Rapoport condition, where $M_{a,b}$
   (resp. $M_{b,a}$) is a \dieu submodule of single slope $\frac{b}{a+b}$
   (resp. $\frac{a}{a+b}$). Then $a$ is a multiple of $e$. Therefore,
   one can not construct a special \dieu module, in the sense of
   Manin, with an 
   $\O$-action which satisfies the Rapoport condition and has arbitrary
   possible slope sequences in Lemma~\ref{32}. The proof of this fact
   is as follows.  
  
  Let $M_{a,b}=M^0\oplus \dots \oplus M^{f-1}$ as in (\ref{22}), where
  $M^i$ 
  is a free $W^i$-module of rank one. We can choose a basis $x_i$ for
  $M^i$ such that $Vx_{i+1}=\pi^{n_i} x_i$ for $1\le i \le f-1$ and
  $Vx_1=u_0\pi^{n_0}x_0$, for some $u_0\in (W^0)^\times$. As
  $M$ satisfies the Rapoport condition, we have $n_i=0$ or $e$. 
  As the slope of $M_{a,b}$ is
  $\frac{b}{a+b}$, we have $V^fx_0=\pi^a u x_0$ for some unit $u$ 
  hence $a=n_0+\dots+n_{f-1}$, which is a multiple of $e$.    
\end{rem}

\subsection{}
Let $V$ be a $2$-dimensional $\FF_v$-vector space with a non-degenerate
alternating form $\psi$ on $V$ with values in $\Q_p$ such that
$\psi(ax,y)=\psi(x,ay)$ for $a\in \FF_v,\ x ,y\in V$. Let $G$ be the
algebraic group of $\FF_v$-linear similitudes over $\Q_p$ and let $G_1$ be the
derived group of $G$. The algebraic group $G_1$ is simply-connected
and it is the kernel of the multiplier map $c:G\to\Gm$.

Let $M$ be a quasi-polarized \dieu $\O$-module. 
We choose an isomorphism between $M\otimes_W B(k)$ and $V\otimes_{\Q_p}
B(k)$ for skew-symmetric $\FF_v$-modules over $B(k)$. Let 
$b\in G(B(k))$ be the element obtained by the transport of structure of the
Frobenius $F$ on $M$. The $\sigma$-conjugacy class $B(G)$ classifies
the ($F$-)isocrystals with $G$-structure. The fibre of $\nu(b)$
under the Newton map $\nu:B(G)\to S(g)$ is classified by $H^1(\Q_p,J_b)$ [RR,
Prop.~1.17], where $J_b$ be the algebraic group over 
$\Q_p$ which represents the group functor [RZ, Prop.~1.12]
\[ R\mapsto \{g\in G(R\otimes_{\Q_p}B(k)); g(b\sigma)=(b\sigma)g
\}. \]

\begin{lemma}\label{35}
 $ H^1(\Q_p, J_b)=0.$
\end{lemma}
\begin{proof} 
Let $N:=M\otimes B(k)$ be the isocrystal with $G$-structure. The group
$J_b$ is $\Aut_G(N)$, viewed as an algebraic group over $\Q_p$. 
If $M$ is supersingular, then $\Aut_G(N)$ is the multiplicative group of
an quaternion algebra over $\FF_v$ with reduced norm in $\Q_p$. 
If $M$ is not supersingular, then $\Aut_G(N)=\FF_v^\times \times
\Q^\times_p$. In both cases, $H^1(\Q_p,J_b)=0$.
\end{proof}

\begin{cor}\label{36}
  Let $M_1$ and $M_2$ are two quasi-polarized \dieu $\O$-modules. If
  $\slope(M_1)=\slope(M_2)$, then $M_1\otimes_W B(k)$ and
  $M_2\otimes_W B(k)$ are isomorphic as quasi-polarized isocrystals
  with the action by $\FF_v$. 
\end{cor}

\begin{cor}\label{38}
  Any polarized abelian $O_\FF$-variety over $k$ 
  is isogenous to one which satisfies the Rapoport condition.
\end{cor}
\begin{proof}
  It follows from (\ref{33}) and Corollary~\ref{36} that the statement
  holds for the associated $p$-divisible group. 
  Then the assertion follows from a theorem of Tate.  \qed
\end{proof}

\section{Normal forms}
\label{sec:04}

\subsection{}
Let $M$ be a non-ordinary separably quasi-polarized \dieu $\O$-module
over $k$ which 
satisfies the Rapoport condition. Let $\ul a(M)=(a^i)_i$ be the $a$-type
of $M$. Note that the Lie type of $M$ is $(\{0,e\})_i$ and the
constraint for the $a$-type is $0\le a^i\le e$. Denote by $\tau(M)$ the \emph{$a$-index} of $M$, which is defined to
be the subset of $\Z/f\Z$:
\[ \tau(M):=\{i\in \Z/f\Z; a^i\not= 0\}. \] 
Write $\tau$ for $\tau(M)$. 

Let $\calD^{-1}=(\pi^{-d})$ be the inverse of the different of $\O$
over $\Z_p$. There is a unique $W\otimes \O$-bilinear pairing
$ (\, , ): M\times M\to W\otimes  \O$
such that $\<x,y\>=\Tr_{W\otimes\O /W} (\pi^{-d}(x,y))$. 
Write $\ol M:=M/\pi M$. The module $\ol
M$ is a $2f$-dimensional vector space over $k$ together with a
$k$-linear action by
$\O/\pi=\O^{\mathrm{ur}}/p$ that commutes with the action of $F$ and
$V$. 
The perfect pairing $(\, , )$ on $M$ induces an perfect pairing
$(\, , )$ on $\ol M$ that satisfies $(Fx,y)=(x,Vy)^p$
and $(ax,y)=(x,ay)$ for all 
$x,y \in \ol M, a\in \O^{\mathrm{ur}}/p$. 
The decomposition of $\ol M$ (\ref{22}) into $\sigma_i$-eigenspaces  
\[ \ol M=\oplus_{i\in \Z/f\Z} \ol M^i \]
respects the perfect pairing as before (\ref{22}). 

As $M$ satisfies the Rapoport condition, $V\ol M$ and $F\ol M$ are
free $k\otimes \O^{\ur}$-module of rank one,  and we have that $F\ol M=\ker
V$ and $V \ol M=\ker F$. Consider $\ol M$ as a $k[F,V]$-module, it is isomorphic to
$N/pN$ for a separably quasi-polarized \dieu  $\O^{\ur}$-module $N$
over $k$ of rank $2f$ with the same $a$-index as $M$. 

\begin{prop}
  There exists a $k$-basis $\{x_i, y_i\}$ of $\ol M^i$ for each $i\in
  \Z/f\Z$ such that
  \begin{itemize}
  \item $(x_i,y_i)=1$ and $y_i\in V\ol M$ for all $i\in \Z/f\Z$,
  \item $F\, x_i=
      \begin{cases}
        x_{i+1} & \text{if $i+1\not\in \tau$}\\
        -y_{i+1}  &
           \text{if $i+1\in \tau$}
      \end{cases}\quad $ 
      $V y_i=
      \begin{cases}
        y_{i-1} & \text{if $i\not\in \tau$}\\
        0 & \text{if $i\in \tau$} 
      \end{cases}\quad $  
      $V x_i=
      \begin{cases}
        y_{i-1} & \text{if $i\in \tau$}\\
        0 & \text{if $i\not\in \tau$} 
      \end{cases}$ 
  \end{itemize}
\end{prop}
\begin{proof}
  It follows from [Y1, Prop. 4.1]. \qed
\end{proof}

Note that this proposition gives the classification of the $\pi$-torsion
subgroup scheme $A[\pi]$ of separably polarized abelian $O_\FF$-varieties
over $k$, classified by the $a$-indices. \\

The following results generalize [Y1, Prop.~4.1, Prop.~4.2, Lemma 4.3],
the proofs are the same and omitted.

\begin{prop}\label{normal}
  There exists a $W^i$-basis $\{X_i,Y_i\}$ of $M$ for each $i\in
  \Z/f\Z$ such
  that
  \begin{itemize}
    \item[(i)] $(X_i, Y_i)=1$ and $Y_i\in (VM)^i,\ \forall \,i\in
      \Z/f\Z$, 
    \item[(ii)] $F\, X_i=
      \begin{cases}
        X_{i+1} & \text{if $i+1\not\in \tau$}\\
        -Y_{i+1}+c_{i+1}\pi X_{i+1}  &
           \text{if $i+1\in \tau$} 
      \end{cases}$\quad\quad% for some $a_{i+1}\in W^{i+1}$, \\
      $F\, Y_i=
      \begin{cases}
        p Y_{i+1} & \text{if $i+1\not\in \tau$}\\
        p X_{i+1} & \text{if $i+1\in \tau$}
      \end{cases}$ \\for some $c_{i+1}\in W^{i+1}$.
  \end{itemize}
\end{prop}

\begin{prop}\label{normal2}
  There exists a
  $W^i$-basis $\{X_i,Y_i\}$ of $M^i$ for each $i\in \Z/f\Z$ such that
  \begin{itemize}
    \item[(i)] $(X_i,Y_i)\in (W^i)^\times$ and $Y_i\in (VM)^i, \
      \forall \, i\in \Z/f\Z$,
    \item[(ii)] $F\, X_i=
      \begin{cases}
        X_{i+1} & \text{if $i+1\not\in \tau$}\\
        Y_{i+1}+c_{i+1}\pi X_{i+1} &
           \text{if $i+1\in \tau$ }
      \end{cases}$  % \\ for some $a_{i+1}\in W(k)$\\
      $F\, Y_i=
      \begin{cases}
        \pi^e Y_{i+1} & \text{if $i+1\not\in \tau$}\\
        \pi^e X_{i+1} & \text{if $i+1\in \tau$,}
      \end{cases}$ \\ for some $c_{i+1}\in W^{i+1}$
  \end{itemize}
\end{prop}

Note that the $i$-th component $a^i$ of $\ul a(M)$
is $\min\{e,\ord_\pi(c_i\pi)\}$ for $i\in\tau$.

\begin{lemma}\label{45}
If $M$ is superspecial (\ref{222}) (in this case $a^i=e$, $\forall\,
i\in \Z/f\Z$) , then 

(1) There exists a $W^i$-basis $\{X_i, Y_i\}$ of $M^i$ for each $i\in
    \Z/f\Z$ such that 
  \begin{itemize}
  \item $Y_i\in (VM)^i$ and $(X_i,Y_i)=1$,
  \item $FX_i=-Y_{i+1},\ FY_i=pX_{i+1}$,
  \end{itemize}
for all $i\in \Z/f\Z$.

(2) There exists a $W^i$-basis $\{X_i,Y_i\}$ of $M^i$ for each
    $i\in\Z/f\Z$ such that
  \begin{itemize}
  \item $Y_i\in (VM)^i$ and $(X_i,Y_i)=u_i\in (W^i)^\times $
    with $u_i^{\sigma^r}=u_i$, where $r=\mathrm{lcm}(2,f)$,
  \item $FX_i=Y_{i+1},\ FY_i=pX_{i+1}$,
  \end{itemize}
for all $i\in \Z/f\Z$.
\end{lemma}

For applications, we need to classify all quasi-polarized superspecial 
\dieu $\O$-modules, not just separably-polarized ones or those
satisfying the Rapoport condition. 

\begin{lemma}\label{455}
  Let $M$ be a quasi-polarized superspecial \dieu $\O$-module over $k$
and let $e_1, e_2$ be as in (\ref{222}). 

(1) If $f$ is even, then there is a $W^i$-basis $X_i, Y_i$ for $M^i$
    for each $\iz$ such that
    \begin{itemize}
    \item [(i)] $(X_i,Y_i)=
      \begin{cases}
        \pi^n & \text{if $i$ is odd},\\
        \pi^{n+e-e_1-e_2} & \text{if $i$ is even}, 
      \end{cases}$ \\ for all $\iz$ and some $n\in \Z$.
    \item [(ii)] $FX_i=
      \begin{cases}
        -\pi^{e_1}Y_{i+1} & \text{if $i$ is odd},\\
        -\pi^{e-e_2}Y_{i+1} & \text{if $i$ is even},\\
      \end{cases}\quad FY_i=
       \begin{cases}
        v\pi^{e_2}X_{i+1} & \text{if $i$ is odd},\\
        v\pi^{e-e_1}X_{i+1} & \text{if $i$ is even},\\
      \end{cases}$ \\ where $v\pi^e=p$.    
    \end{itemize}

(2) If $f$ is odd, then there is a $W^i$-basis $X_i, Y_i$ for $M^i$
    for each $\iz$ such that
    \begin{itemize}
    \item [(i)] $(X_i,Y_i)=\pi^n$ for some $n\in \Z$.
    \item [(ii)] $FX_i=-\pi^{e_1}Y_{i+1}, FY_i=
        v\pi^{e_2}Y_{i+1}$ for $\iz$, where $v\pi^e=p$.
    \end{itemize}
\end{lemma}
\begin{proof}
  The proof is similar to that of [Y1, Lemma 4.3], hence is sketched. 

  (1) Write $f=2c$ and let ${M'}:=\{x\in M\, |\, F^c  
  x=(-1)^c V^cx\, \}$. Since $M$ is superspecial, we have $F^2
  {M'}={pM'}$ and ${M'}\otimes_{W(\F_{p^f})} W(k)\simeq M$. We can
  choose bases $\{X_0,Y_0\}$, $\{X_1, Y_1\}$ for $(M')^0, (M')^1$
  respectively such that $FX_0=-\pi^{e-e_2}Y_{1},FY_0=
  v\pi^{e-e_1}X_{1}$, and $(X_1,Y_1)=\pi^n$ for some $n\in \Z$. 

  Define $X_i,Y_i$ recursively for $2\le i\le f$ by (ii). Then it is
  straight forward to verify that $X_f=X_0$, $Y_f=Y_0$ and (i).

  (2) Write $f=2c+1$ and let $M':=\{
  x\in M\, |\, F^{2f}x+p^f x=0\}$. From (ii), $Y_0$ is required to be
  $(-1)^{c+1}\pi^{-e_1} p^{-c}F^fX_0$. Let $X_0,Z_0$ be a basis for
  $(M')^0$ such that $FX_0=\pi^{e_1}X'_1, FZ_0=\pi^{e_2}Y'_1$ for some
  basis $\{X'_1,Y'_1\}$. Let $Y_0:=(-1)^{c+1}\pi^{-e_1} p^{-c}F^fX_0$. 
  If $e_1<e_2$, then $X_0, Y_0$ form a basis. If $e_1=e_2$, then we can
  choose again $X_0$ such that $X_0, Y_0$ form a basis for
  $(M')^0$. Define $X_i,Y_i$ recursively for $1\le i\le f$ by (ii)
  and it is easy to check that $X_f=X_0$ and $Y_f=Y_0$. We found a
  basis satisfying (ii).

  Write $(X_0,Y_0)=\mu \pi^n$ for some unit $\mu$ and some integer
  $n$. It follows from $(F^fX_0,F^fY_0)=p^f(X_0,Y_0)^{\sigma^f}$ that
  $\mu^{\sigma^f}= \mu$.      
  If we replace $X_0$ by $\lambda X_0$, then
  $\mu$ will change to $\lambda \lambda^{\sigma^f}\mu$. Since
  $B(\F_{p^{2f}})[\pi]/B(\F_{p^f})[\pi]$ is a quadratic unramified extension,
  we can adjust $X_0$ by choosing a suitable $\lambda$ such that
  $(X_0,Y_0)=\pi^n$, hence (i) is satisfied. \qed
\end{proof}
\subsection{}
\label{46}
Let $M$ be a non-ordinary separably quasi-polarized \dieu $\O$-module
satisfying the Rapoport condition. Let $\ul a(M)=(a^i)$ be the $a$-type
of $M$ and 
$\tau=\tau(M)=\{n_1, n_2,\dots n_t\}$ be the $a$-index of $M$, where
$t=|\ul a(M/\pi M)|$, the \emph{reduced $a$-number of $M$}. We
assume that $0\le n_1< n_2<\dots<n_t<f$ and let
$n_0:=n_t-f$ and $n_{t+1}=n_1$. Denote by 
$A(e,f)$ the set of possible $a$-types on the Rapoport locus
\[ A(e,f):=\{(a^i); i\in \Z/f\Z, a^i\in \Z, 0\le a^i\le e \}, \]
with the partial order that $(a^i)\le (b^i)$ if $a^i\le b^i,\
\forall\, i\in \Z/f\Z$.  
An $a$-type $(a^i)\in A(e,f)$ is called \emph{spaced} if $a^i
a^{i+1}=0, \ \forall i\in \Z/f\Z$, cf. [GO, Sect. 1, p.112]. 
 
By Proposition~\ref{normal2}, we can choose a $W^i$-basis $\{X_i,
Y_i\}$ of $M^i$ for each $i\in \Z/f\Z$ such that
\[ 
\begin{pmatrix}
  F X_{i-1} \\ FY_{i-1} 
\end{pmatrix}=A_{i}
\begin{pmatrix}
  X_{i} \\ Y_{i}
\end{pmatrix}, \]
where 
\[ A_{i}=
\begin{pmatrix}
  1 & 0 \\ 0 & \pi^e 
\end{pmatrix}
\ \text{if $i\not \in \tau$}; \quad A_{i}=  
\begin{pmatrix}
  \pi^{a^{i}}c_i & 1 \\
  \pi^e & 0 
\end{pmatrix}\ \text{if $i\in \tau$},   \]
for some $c_i\in W^i$ for $i\in \tau$ and $c_i$ is a unit if $a^i<e$. 

Let $\ell_i:=n_i-n_{i-1}$ for $1\le i\le t$. We have
\[ 
\begin{pmatrix}
  F^{\ell_{i}} X_{n_{i-1}} \\ F^{\ell_{i}} Y_{n_{i-1}}
\end{pmatrix}=
\begin{pmatrix}
    \pi^{a^{i}}c_{n_{i}} & 1 \\
  \pi^{e\ell_{i}} & 0 
\end{pmatrix}
\begin{pmatrix}
  X_{n_{i}} \\ Y_{n_{i}}
\end{pmatrix},\quad \forall 1\le i\le t. \]

\nc{\xs}{X_{n_{s-1}}}
\nc{\xsm}{X_{n_{s-2}}}
\nc{\fs}{F^{\ell_s+\ell_{s+1}+\dots+\ell_t}}
\nc{\fsm}{F^{\ell_{s-1}+\ell_s+\dots+\ell_t}}
\nc{\ys}{Y_{n_{s-1}}}
\nc{\ysm}{Y_{n_{s-2}}}
\nc{\xone}{X_{n_t}}
\nc{\yone}{Y_{n_t}}

%\nc{\asp}{\alpha_s'}
%\nc{\bsp}{\beta_s'}
%\nc{\csp}{\gamma_s'}
%\nd{\dsp}{\delta_s'}

For $1\le s\le t$, suppose that 

\[  \begin{pmatrix}
    \fs \xs \\\fs \ys
  \end{pmatrix}=
  \begin{pmatrix}
    \alpha_s & \beta_s  \\
    \gamma_s & \delta_s
  \end{pmatrix}
  \begin{pmatrix}
    \xone \\ \yone 
  \end{pmatrix},
\]
for some coefficients $\alpha_s, \beta_s,\gamma_s, \delta_s$. It
follows from 
\begin{equation}
  \label{eq:46}
  \begin{split}
    \begin{pmatrix}
  \fsm \xsm \\ \fsm \ysm 
\end{pmatrix} & =
\begin{pmatrix}
  \fs (F^{\ell_{s-1}} \xsm) \\ \fs (F^{\ell_{s-1}}\ysm)
\end{pmatrix} \\
&= \begin{pmatrix}
    \pi^{a^{n_{s-1}}}c_{n_{s-1}}^{(\ell_s+\ell_{s+1}+\dots+\ell_t)} & 1 \\
  \pi^{e\ell_{s-1}} & 0 
\end{pmatrix}
\begin{pmatrix}
  \fs \xs \\ \fs \ys 
\end{pmatrix}
  \end{split}
\end{equation}
that 
\begin{equation}
  \label{eq:47}
  \begin{pmatrix}
    \alpha_{s-1} & \beta_{s-1}  \\
    \gamma_{s-1} & \delta_{s-1}
  \end{pmatrix}=
  \begin{pmatrix}
    u_{s-1} & 1 \\
  \pi^{e\ell_{s-1}} & 0
  \end{pmatrix} 
  \begin{pmatrix}
    \alpha_s & \beta_s  \\
    \gamma_s & \delta_s
  \end{pmatrix},\quad    
   \begin{pmatrix}
    \alpha_{t} & \beta_{t}  \\
    \gamma_{t} & \delta_{t}
  \end{pmatrix}=
\begin{pmatrix}
    \pi^{a^{n_t}}c_{n_t} & 1 \\
  \pi^{e\ell_{t}} & 0 
\end{pmatrix},
\end{equation}
where
$u_{s-1}=\pi^{a^{n_{s-1}}}c_{n_{s-1}}^{(\ell_s+\ell_{s+1}+\dots+\ell_t)}$.
Recall that we write $a^{(n)}$ for $a^{\sigma^n}$ (\ref{33}). Therefore, 
\begin{equation}
  \label{eq:461}
  \begin{pmatrix}
    \alpha & \beta  \\
    \gamma & \delta
  \end{pmatrix}:=   
  \begin{pmatrix}
    \alpha_1 & \beta_1  \\
    \gamma_1 & \delta_1
  \end{pmatrix}=   
  \begin{pmatrix}
    u_{1} & 1 \\
    \pi^{e\ell_{1}} & 0
  \end{pmatrix}     
  \begin{pmatrix}
    u_{2} & 1 \\
    \pi^{e\ell_{2}} & 0
  \end{pmatrix}\cdots
  \begin{pmatrix}
    u_{t} & 1 \\
  \pi^{e\ell_{t}} & 0
  \end{pmatrix}.
\end{equation}

\subsection{}
\label{47}
Consider the case that $c_i=0$ for all  $i\in \tau$. If $t$ is even,
write $t=2d$, then we have 
\[    \begin{pmatrix}
    \alpha & \beta  \\
    \gamma & \delta
  \end{pmatrix}=
  \begin{pmatrix}
    \pi ^{e(\ell_2+\ell_4+\dots+\ell_{2d})} & 0 \\
    0 & \pi ^{e(\ell_1+\ell_3+\dots+\ell_{2d-1})}
  \end{pmatrix} \]
and $\slope(M)=s(i)$ (\ref{321}), where $i=\min\{e(\sum_{1\le j\le d}
\ell_{2j}), 
e(\sum_{1\le j\le d}\ell_{2j-1})\}$. If $t$ is odd, write $t=2d+1$,
then we have
\[  \begin{pmatrix}
    \alpha & \beta  \\
    \gamma & \delta
  \end{pmatrix}^2=
  \begin{pmatrix}
    \pi ^g & 0 \\
    0 & \pi^g  
  \end{pmatrix} \] 
and $M$ is supersingular.
\begin{prop}\label{48}
  Let $M$ be a separably quasi-polarized \dieu $\O$-module over
  $k$ that satisfies the Rapoport condition. If $\ul a(M)$ is spaced,
  i.e. $\ell_i\ge 2\ \forall\, 1\le i\le t$, 
  then $\slope(M)\ge s(|\ul a(M)|)$ (\ref{321}). 
\end{prop}
\begin{proof}
  We may assume that $M$ is non-ordinary. It follows from
  (\ref{eq:461}) that 
  \[ \begin{pmatrix}
  F^f \pi^e \xone \\ F^f \yone 
  \end{pmatrix}=A \begin{pmatrix}
   \pi^e \xone \\ \yone 
  \end{pmatrix}, \quad A=\prod_{i=1}^t A_i= \pi^{|\ul a(M)|} \prod_{i=1}^t
  \left(\pi^{-a^{n_i}}A_i\right), \quad A_i= 
  \begin{pmatrix}
    u_{i} & \pi^e \\
  \pi^{e(\ell_{i}-1)} & 0
  \end{pmatrix}. \]
Let $N$ be the \dieu $\O$-submodule of $M$ generated by $\pi^e \xone,
\yone$. As $F^f (N)\subset \pi^{|\ul a(M)|}N$, we get $\slope(M)\ge
s(|\ul a(M)|)$. \qed
\end{proof}

In general the slope sequence $\slope(M)$ is not determined by $\ul
a(M)$. This is known even in the unramified case [GO]. We will see in
(\ref{68}) that the bound in (\ref{48}) is sharp for spaced $a$-types,
while it is not the case for non-spaced ones, see (\ref{610}).
\subsection{}\label{49}
When $t=|\tau|=1$, say $\tau=\{0\}$, we have
\[ 
\begin{pmatrix}
  F^f X_0 \\ F^f Y_0 
\end{pmatrix}=
\begin{pmatrix}
  c_0 & 1 \\ \pi^g & 0 
\end{pmatrix}
\begin{pmatrix}
  X_0 \\ Y_0 
\end{pmatrix}.
\]
for some $c_0\in \pi^{a^0} W^0$. 
It is easy to see that $M$ satisfies a Cayley-Hamilton equation ([O2],[Y1])
$F^{2f}X_0-c_0^{(f)}F^fX_0-\pi^g X_0=0$, thus $\slope(M)=s(i)$, 
where $i=\min\{\frac{g}{2},\ord_\pi(c_0)\}$.

\subsection{}
\label{410}
When $t=|\tau|=2$, say $\tau=\{0,\ell_2\}$ and $\ell_2\le \ell_1$, we
have 
\nc{\xl}{X_{\ell_1}}
\nc{\yl}{Y_{\ell_1}}
\[   \begin{pmatrix}
   F^f \xl \\ F^f \yl  
  \end{pmatrix}=
  \begin{pmatrix}
    u_1 & 1 \\ \pi^{e\ell_1} & 0  
  \end{pmatrix}
  \begin{pmatrix}
    u_2 & 1 \\ \pi^{e\ell_2} & 0  
  \end{pmatrix}
  \begin{pmatrix}
    \xl \\ \yl
  \end{pmatrix}, \]
for some coefficients $u_1, u_2$ ($u_1=c_0^{(\ell_2)},
u_2=c_{\ell_2}$). Then
\begin{equation}
  \label{eq:48}
  \begin{pmatrix}
    \alpha & \beta \\ \gamma & \delta
  \end{pmatrix}=
   \begin{pmatrix}
    u_1 u_2+\pi^{e\ell_2} & u_1 \\
    \pi^{e \ell_2}u_2  & \pi^{e \ell_2} \\
  \end{pmatrix}
\end{equation}

If $u_1=0$, then the matrix is $
\begin{pmatrix}
    \pi^{e\ell_2} & 0 \\
    \pi^{e \ell_1}u_2  & \pi^{e \ell_1} \\
\end{pmatrix}$, hence $\slope(M)=s(e\ell_2)$. If $u_1\not=0$, then 
    we have the Cayley-Hamilton equation
\[ \beta F^{2f}
\xl-(\beta^{(f)}\delta+\beta\alpha^{(f)})F^f\xl+\beta^{(f)}(\alpha\delta-\beta
\delta)\xl=0. \]
As $\ord_\pi(\delta)=e\ell_1\ge \frac{g}{2}$, we see that
$\slope(M)=s(i)$, where $i=\min\{ \frac{g}{2},\ord_\pi(a^{(f)}) \}$.
In both cases we have that $\slope(M)=s(i)$, where $i=\min\{
\frac{g}{2},\ord_\pi(u_1u_2+\pi^{e\ell_2}) \}$. 

\section{Alpha stratification}
\label{sec:05}

\subsection{}\label{51}
Let $p$ be a fixed prime number. Let $\scrM^{\DP}$ denote the moduli
stack over $\Spec \Z_{(p)}$ of polarized abelian $O_\FF$-varieties
$(A,\lambda,\iota)$ of dimension $g=[F:\Q]$ with the polarization
$\lambda$ of prime-to-$p$ degree. It is a separated Deligne-Mumford
algebraic stack over $\Spec \Z_{(p)}$ locally of finite type. In [DP],
Deligne and Pappas showed that the
algebraic stack $\scrM^{\DP}$ is flat and a locally 
complete intersection over $\Spec \Z_{(p)}$ of relative dimension $g$,
and the closed fibre $\scrM^{\DP}\otimes \Fp$ is geometrically normal
and has singularities of codimension at least two. It follows from
Deligne-Pappas' results and the compactification of Rapoport that the
irreducible components of geometric 
special fibre $\scrM^{\DP}\otimes \ol \Fp$ are in bijection correspondence with
those of geometric generic fibre $\scrM^{\DP}\otimes \Qbar$. Those are
parameterized by the isomorphism classes of non-degenerate skew-symmetric
$O_\FF$-modules $H_1(A(\C), \Z)$ for all $(A,\lambda,\iota)\in
\scrM^{\DP}(\C)$.   

Let $\scrM^R$ denote the Rapoport locus of $\scrM^{\DP}$. It is the
smooth locus of $\scrM^{\DP}$ (\ref{215}). Let $\scrM$ denote the
reduction $\scrM^R\otimes_\Z k(v)$ of $\scrM^R$ modulo $v$, where
$k(v)$ is the residue field of $O_\FF$ at $v$. We will define the
stratification on $\scrM$ by $a$-types scheme-theoretically, (cf.
[Y1, Sect. 3]).

\subsection{}\label{52}
Let $S$ be a locally noetherian scheme over $\Spec k(v)$ and
$\pi_A:(A,\lambda,\iota)\to S\in \scrM(S)$ a polarized abelian $O_\FF$-scheme
over $S$. The sheaf $R^1(\pi_A)_*(\O_A)$ is a locally free rank one $O_\FF\otimes
\O_S$-module. It admits a decomposition (\ref{22})
\[ R^1(\pi_A)_*(\O_A)=\bigoplus_{i\in \Z/f\Z} R^1(\pi_A)_*(\O_A)^i \]
with respect to the action by $O_\FF$. Each component
$R^1(\pi_A)_*(\O_A)^i$ is a locally free $O_S[\pi]/(\pi^e)$-module of rank
one and it admits a filtration of locally free $O_S$-modules
\[ R^1(\pi_A)_*(\O_A)^i=\scrF^{i,0}\supset \scrF^{i,1}\supset \dots
\supset \scrF^{i,e}=0, \quad \scrF^{i,j}:=\ker \pi^{e-j}.  \]
Let $F_{A/S}: A\to A^{(p)}$ be relative Frobenius morphism over $S$,
where $A^{(p)}:=A\times_{S, F_{\mathrm{abs}}} S$ and
$F_{\mathrm{abs}}$ is the absolute Frobenius morphism on $S$. It
induces an $\O_S\otimes O_\FF$-linear morphism $F^*_{A/S}:
R^1(\pi_A)_*(\O_A)^{(p)}\to  R^1(\pi_A)_*(\O_A)$. 
  
Let $\ul a=(a^i)\in A(e,f)$ (\ref{46}) be an $a$-type. Let
$\scrM_{\ge \ul a}$ denote the substack of $\scrM$ whose objects
$\pi_A:(A,\lambda, \iota)\to S$ satisfy the following condition
\begin{equation}
  \label{eq:52}
  F^*_{A/S}\left(R^1(\pi_A)_*(\O_A)^{(p),i}\right)\subset
  \scrF^{i,a^i}, \quad \forall\, i\in \Z/f\Z.
\end{equation}
It is clear that the condition (\ref{eq:52}) is closed and locally for
Zariski topology $\scrM_{\ge \ul a}$ is defined by finitely many
equations. Therefore $\scrM_{\ge \ul a}$ is a closed algebraic
substack of $\scrM$. 

It is clear that $\scrM_{\ge \ul a}$ is non-empty
as it contains $E\otimes O_\FF$ with a prime-to-$p$ polarization, where
$E$ is a supersingular elliptic curve.

\begin{lemma}\label{53}
  Let $(A,\lambda,\iota)\in \scrM(k)$ and $\ul a(A)=(a^i)$ be the $a$-type of
  $A$ (\ref{27}). Then $\coker F^*_{A/k}=\oplus_{\iz}
  k[\pi]/(\pi^{a^i})$. 
\end{lemma}
\begin{proof}
  Let $M^*$ denote the contravariant \dieu module of $A$. We identify
$M^*/p M^*=H^1_{DP}(A/k)$ and have $H^1(A,\O_A)=M^*/VM^*$ and $\coker
F^*_{A/k}=M^*/(F,V)M^*$. On the other hand, the quasi-polarization
gives an isomorphism $M\simeq M^*$ of \dieu $\O$-modules, which
induces an isomorphism $M/(F,V)M\simeq M^*/(F,V)M^*$ of 
$O_\FF\otimes k$-modules. \qed 
\end{proof}

\begin{thm}\label{54}
  The algebraic stack $\scrM_{\ge \ul a}$ is smooth over $\Spec k(v)$
  of pure dimension $g-|\,\ul a\,|$.
\end{thm}
\begin{proof}
Let $x:\Spec k\to \scrM_{\ge \ul a}$ be a geometric point.
  Then there is an affine open neighborhood $U$ of $x$ and a polarized
  abelian $O_\FF$-scheme $(A,\lambda,\iota)\in \scrM(U)$ 
  such that $R^1(\pi_A)_*(\O_A)^{(p)}$ and $R^1(\pi_A)_*(\O_A)$ are
  free $\O_U\otimes O_\FF$-modules. Let $x_i^{(p)}$ and $x_i$ be
  $O_U[\pi]/(\pi^e)$-bases of $R^1(\pi_A)_*(\O_A)^{(p),i}$ and
  $R^1(\pi_A)_*(\O_A)^i$ respectively for each $i\in\Z/f\Z$. Let
  \[ F^*_{A/S}(x^{(p)}_i)=\sum_{0\le j < e} f_{i,j} \pi^j x_i, \quad
  f_{i,j}\in \O_U.\] 
  Then locally
  $\scrM_{\ge \ul a}$ is defined by the equations $f_{i,j}$,
  where $(i,j)\in I=\{(i,j); i\in \Z/f\Z, 0\le j < a^i\}$. Therefore,
  $\dim_x \scrM_{\ge \ul a}\ge g-| I|=g-|\,\ul a\,|$.

Let $\ol M:=H^1_{cris}(A_x/k),\ol {M}^{(p)}:=H^1_{cris}(A_x^{(p)}/k)
=\ol M\otimes_{k,\sigma} k$. We have the Hodge filtration:
\[ \begin{CD}
0 @>>> Fil^{(p)} @>>> \ol{M}^{(p)} @>>>
  Q^{(p)}:=H^1(A_x,\O_{A_x})^{(p)} @>>> 0\\ 
  @.             @. @VV{F^*}V @VV{F^*}V \\
0 @>>> Fil @>>> \ol M @>>> Q:=H^1(A_x,\O_{A_x}) @>>> 0\\
\end{CD}\]
By Proposition~\ref{normal} and Lemma~\ref{53}, we can choose a
$k[\pi]/(\pi^e)$-basis $\{x_i,y_i\}$ be of $(\ol M)^i$ for 
each $i\in \Z/f\Z$ such that  $y_i\in Fil^i$ and $F^*(x_i\otimes
1)=-y_i+c_i x_i$ for some $c_i\in \pi^{a^i}k[\pi]/(\pi^e)$. 
For any PD-extension $R$ of $k$, the differential forms $x_i,y_i$ have
unique horizontal  
liftings in $\wt{ M}:=H^1_{cris}(A/R)$ with respect to the
Gauss-Manin connection, which we will denote the liftings by $x_i,y_i$ again. 
\ncmd{\sspan}{\mathrm{span}} 

Let $R:=k[[\underline{t}]]/(\underline t)^2$ be the first order
universal deformation ring of $(A_x,\lambda_x,\iota_x)$, where $\ul
t=(t_{i,j}), i\in \Z/f\Z, 0\le j<e$. The first order universal
deformation $(\tilde{A},\tilde{\lambda},\tilde{\iota})$ gives rise to
 \[ \begin{CD}
0 @>>> \widetilde{Fil}^{(p)} @>>> \wt{M}^{(p)} @>>> 
\wt{Q}^{(p)} @>>> 0\\
  @.             @. @VV{F^*}V @VV{F^*}V \\
0 @>>> \widetilde{Fil} @>>> \wt{M} @>>> \wt{Q} @>>> 0\\
\end{CD}\] 
where $\widetilde{Fil}=\sspan\<y_i+\sum_{0\le j<e}t_{i,j}\pi^j x_i\>, \ 
\widetilde{Fil}^{(p)}=\sspan\<(y_i+\sum_{0\le j<e}t_{i,j}\pi^j x_i)\otimes 1\>
=\sspan\<y_i\otimes 1\>, \ 
\wt{Q}=\sspan\<[x_i]\>$, and $\wt{Q}^{(p)}=\sspan\< [x_i\otimes
1]\>$ for all $\iz$. Here the bracket $[\ ]$ denotes the class modulo
$\wt{Fil}$ and $\wt{Fil}^{(p)}$. It follows from 
\[ F^*([x_i\otimes
1])=[-y_i]=\sum_{0\le j<e} t_{i,j}[\pi^j x_i]\] 
that $R/(t_{i,j})_{(i,j)\in I}$ is the first order deformation ring of   
$\scrM_{\ge \ul a}$ at $x$ and the tangent space has dimension
$g-|\,\ul a\, |\le \dim_x\scrM_{\ge \ul a}$. The assertion follows. \qed
\end{proof}

\begin{cor}
  The ordinary points are dense in
  $\scrM^{\DP}\otimes \F_p$. 
\end{cor}
\begin{proof}
  Note that for ordinary points $a^i=0,\forall\, i\in \Z/f\Z$. The
  statement follows from Theorem~\ref{54} and the density of the Rapoport
  locus.\qed 
\end{proof}

\subsection{}
\label{56}

Let $\scrM_{\ul a}$ denote the subset of $\scrM$ that consists of points with
$a$-type ${\ul a}$. It is a locally closed subset of $\scrM$, hence
regarded as a locally closed algebraic substack of $\scrM$ with the reduced
induced structure. Lemma~\ref{53} says that
$(A,\lambda,\iota)\in \scrM_{\ge \ul a}(k)$ if and only if $\ul a(A)\ge
\ul a$. It follows from Theorem~\ref{54} that $\scrM_{\ul a}$ is a
dense open substack of $\scrM_{\ge \ul a}$ and $\scrM_{\ge \ul a}$ is
the scheme-theoretic closure of $\scrM_{\ul a}$ in $\scrM$. This
justifies our notation. 

\section{Deformations of \dieu modules}
\label{sec:06}

\subsection{}
\ncmd{\Cart}{\mathrm{Cart}_p}
We follow the convenient setting of [N, Sect.~0] and [CN, Sect.~2,
p.~1011]. As we will only deal with smooth functors, the deformation
theory developed by P. Norman [N] and Norman-Oort [NO] is enough for
our purpose. We refer the reader to [Z2] for the generalized theory of
displays over more general base ring.  

Let $R$ be a commutative ring of \ch $p$. Let $W(R)$ denote the ring of
Witt vectors over $R$, equipped with the Verschiebung $\tau$ and
Frobenius $\sigma$:
\begin{equation*}
  \begin{split}
    (a_0,a_1,\dots)^\tau&=(0, a_0, a_1,\dots) \\
    (a_0,a_1,\dots)^\sigma&=(a_0^p, a_1^p,\dots).
  \end{split}
\end{equation*}
Let $\Cart(R)$ denote the Cartier ring $W(R)[F][\![V]\!]$ modulo the
relations
\begin{itemize}
\item $FV=p$ and $VaF=a^\tau$,
\item $Fa=a^\sigma F$ and $Va^\tau=aV, \ \forall\, a\in W(R).$
\end{itemize}

A left $\Cart(R)$-module is \emph{uniform} if it is complete and
  separated in the $V$-adic topology. A uniform $\Cart(R)$-module $M$ is
  \emph{reduced} if $V$ is injective on $M$ and $M/VM$ is a free
  $R$-module. A \emph{\dieu module over $R$} is a finitely generated
  reduced uniform $\Cart(R)$-module.

There is an equivalence of categories between the category of finite
dimensional commutative formal group over $R$ and the category of
\dieu module over $R$. We denote this functor by $\bfD_*$. The tangent
space of a formal group $G$ is canonically isomorphic to
$\bfD_*(G)/V\bfD_*(G)$. 

\subsection{}\label{62}

Let $\ul a=(a^i)$ be an $a$-type and $(A_0,\lambda_0,\iota_0)\in
\scrM_{\ge \ul a}(k)$ (\ref{112}) (\ref{52}) be a \emph{non-ordinary} polarized abelian
$O_\FF$-variety. By Lemma~\ref{32}, the associated $p$-divisible group
$G_0=A_0[p^\infty]$ is connected, hence it is a smooth formal group. As
the forgetful functor
$\Def[A_0,\lambda_0,\iota_0]\to\Def[A_0,\iota_0]$ induces an
equivalence of deformation functors, we will consider deformations of
abelian $O_\FF$-varieties and their associated formal groups.  

Let $I$ be the set $\{(i,j); \iz, a^i \le j<e-1
\}$ and let $R:=k[\![t_{i,j}]\!]_{(i,j)\in I}$. We have
\[ \O \otimes_{\Z_p} W(R)= \bigoplus_{\iz} W(R)[T]
/(\sigma_i(P(T))). \]
Set $W^i_{R}:=W(R)[T]/(\sigma_i(P(T)))$ and denote again by
  $\pi$ the image of $T$ in $W^i_{R}$.

Let $M_0$ be the covariant \dieu module of $A_0$. Let $\ul b=(b^i)$ be
the $a$-type of $A_0$ and $\tau$ be the
$a$-index of $M_0$. 
By
Proposition~\ref{normal2}, we can choose a $W^i$-basis of $M^i_0$ for each
$\iz$ such that 
\begin{equation}
  \label{eq:620}
    F X_{i-1}=
      \begin{cases}
        X_{i} & \text{if $i\not\in \tau$}\\
        c_{i} X_{i}+Y_{i} &
           \text{if $i\in \tau$ }
      \end{cases} \quad
      F Y_{i-1}=
      \begin{cases}
        \pi^e Y_{i} & \text{if $i\not\in \tau$}\\
        \pi^e X_{i} & \text{if $i\in \tau$,}
      \end{cases} 
\end{equation}
for some $c_{i}\in \pi^{b^{i}}W^{i}$.

By [NO, Lemma 0.2], we construct a \dieu module $M_{R}$ over
$R$ with a
$\Cart(R)$-linear action by $\O$ as follows. It is the
$\Cart(R)$-module generated by  
$ \{\pi^j X_i, \pi^j Y_i\}_{\iz, 0\le j<e}$ 
with the relations
\begin{align}\label{eq:621}
  F\, X_{i-1}&=
      \begin{cases}
        X_{i} & \text{if $i\not\in \tau$}\\
        c_{i} X_{i}+(Y_{i}+\sum_{a^{i}\le j<e}
        T_{i,j}\pi^j X_{i}) &
           \text{if $i\in \tau$ }
      \end{cases} \\
 Y_{i-1}&=
      \begin{cases}
        V\, u(Y_{i}+\sum_{a^{i}\le j<e}
        T_{i,j}\pi^j X_{i}) & \text{if $i\not\in \tau$}\\
        V\,u\, X_{i} & \text{if $i\in \tau$,}
      \end{cases} \label{eq:622}
\end{align}
where $T_{i,j}$ is the Teichm{\"u}ller lift of $t_{i,j},\ \forall\,
(i,j)\in I$,   
$u=p^{-1}\pi^e$. We impose the natural $W_{R}^i$-module structure on the free
$W(R)$-submodule generated by $ \{\pi^j X_i, \pi^j Y_i\}_{0\le j<e}$,
and impose the trivial $W_{R}^k$-module structure on this submodule for
$k\not=i$. The action of $\O$ on $M_R$ comes from the natural embedding
$\O\hookrightarrow \O\otimes_{\Z_p} W(R)$.

\begin{lemma}
  The \dieu $\O$-module $M_R$ over $R$ is isomorphic to the universal
  deformation of $M_0$ for $\scrM_{\ge \ul a}$.
\end{lemma}
\begin{proof}  
It is clear that $\Cart(k)\otimes_{\Cart (R)}M_R=M_0$. On $M_R/V M_R$,
we have  
\[ F[X_{i-1}]=
        c_{i} [X_{i}]+\sum_{a^{i}\le j<e}
        t_{i,j}\pi^j [X_{i}],\quad \forall \, i\in \tau. \] 
By the Serre-Tate theorem, we obtain a morphism $\Spf R\to 
\scrM_{\ge \ul a,x}^{\wedge}$. By [N, Thm. 1], this construction induces an
        injection of tangent spaces. By Theorem~\ref{54}, 
$\scrM_{\ge \ul a}$ is smooth and $\dim \scrM_{\ge \ul a}=\dim R$. 
Hence the morphism is an isomorphism. \qed
\end{proof}
\begin{rem}
  Another construction of $M_R$ is using tensor products modulo relations. Let
  $P_{R}$ is a free $\O\otimes_{\Z_p} W(R)$-module of rank two,  with
  a $W_{R}^i$-basis $\{X_i, Y_i\}$ for each component $P_{R}^i$. Then
  we construct
  $M_{R}$ to be the quotient of $\Cart(R)\otimes_{W(R)} P_{R}$ 
  modulo the relations (\ref{eq:621}) and (\ref{eq:622}).
\end{rem}

\subsection{}
\label{65}
Let $M_R$ be as in (\ref{62}). Let
\begin{equation}
  \label{eq:651}
   T_i:= 
      \begin{cases}
        \sum_{a^{i}\le j<e}
        T_{i,j}\pi^j & \text{if $i\not\in \tau$}\\
        c_{i}+\sum_{a^{i}\le j<e}
        T_{i,j}\pi^j  &
           \text{if $i\in \tau$.}
      \end{cases} 
\end{equation}
Then we have
\[ 
\begin{pmatrix}
  F X_{i-1} \\ FY_{i-1} 
\end{pmatrix}=A_{i}
\begin{pmatrix}
  X_{i} \\ Y_{i}
\end{pmatrix}, \]
where 
\[ A_{i}=
\begin{pmatrix}
  1 & 0 \\ T_i \pi^e & \pi^e 
\end{pmatrix}
\ \text{if $i\not \in \tau$}; \quad A_{i}=  
\begin{pmatrix}
  T_i & 1 \\
  \pi^e & 0 
\end{pmatrix}\ \text{if $i\in \tau$}.   \]

\begin{lemma}
  The non-ordinary locus of $\scrM_x^{\wedge}$ is defined by
  $\prod_{i\in \tau} t_{i,0}=0$.
\end{lemma}
\begin{proof}
  Take $\ul a=(0,0,\dots, 0)$. On $M_R/(VM_R+\pi M_R)$, we have
\[ F[X_{i-1}]=[X_i]\quad \text{if $i\not \in \tau$}; \quad
F[X_{i-1}]=t_{i,0}[X_i]\quad \text{if $i\in \tau$}. \]
As each slope stratum is reduced, it induces a closed reduced
subscheme,  if not empty, of $\scrM_x^{\wedge}$. 
A point $\grp\in \Spec R$ is in the non-ordinary locus if and only if any of
$t_{i,0}, i\in \tau$, vanishes on $R/\grp$. 
Therefore the defining equation is
$\prod_{i\in \tau} t_{i,0}=0$. \qed
\end{proof}

\subsection{}
\label{67}
We continue with (\ref{65}). Let  
$t=|\ul a(M_0/\pi M_0) |$ and $\tau=\{n_1, n_2,\dots n_t\}$. We
assume that $0\le n_1< n_2<\dots<n_t<f$ and let
$n_0:=n_t-f$ and $n_{t+1}=n_1$. Let $\ell_i:=n_i-n_{i-1}$ for $1\le
i\le t$ (\ref{46}). 

From (\ref{65}), we have 
\begin{equation}
  \label{eq:671}
  \begin{split}
  \begin{pmatrix}
   F^{\ell_{i}} X_{n_{i-1}} \\ F^{\ell_{i}} Y_{n_{i-1}}
   \end{pmatrix}& =
   \begin{pmatrix}
    1 & 0 \\
  T_{n_{i-1}+1}^{(\ell_i-1)} \pi^e & \pi^e 
    \end{pmatrix}
   \dots
   \begin{pmatrix}
      1 & 0 \\
      T_{n_{i-1}+\ell_i-1}^{(1)} \pi^e & \pi^e 
   \end{pmatrix}
   \begin{pmatrix}
      T_{n_i} & 1 \\
      \pi^e & 0
   \end{pmatrix}
   \begin{pmatrix}
      X_{n_{i}} \\ Y_{n_{i}}
   \end{pmatrix} \\
     &=
   \begin{pmatrix}
      1 & 0 \\
      V_i'  & \pi^{e(\ell_i-1)} 
   \end{pmatrix}
   \begin{pmatrix}
      U'_i & 1 \\
      \pi^e & 0
   \end{pmatrix}
   \begin{pmatrix}
     X_{n_{i}} \\ Y_{n_{i}}
   \end{pmatrix} \\
   &= 
   \begin{pmatrix}
     U'_i & 1 \\
     U_i' V_i' +\pi^{e\ell_i} & V_i' 
   \end{pmatrix}
   \begin{pmatrix}
      X_{n_{i}} \\ Y_{n_{i}}
   \end{pmatrix} 
  \end{split}
\end{equation}
where 
\[ U'_i:=T_{n_i},  \quad \text{and}\ \ V'_i:=\sum_{j=1}^{\ell_i-1}
T^{(\ell_i-j)}_{n_{i-1}+j} \pi ^{ej}. \]
Recall that we write $T^{(n)}$ for $T^{\sigma^n}$ (\ref{33}) (\ref{46}).
For $1\le s\le t$, suppose that 

\[  \begin{pmatrix}
    \fs \xs \\\fs \ys
  \end{pmatrix}=
  \begin{pmatrix}
    \alpha_s & \beta_s  \\
    \gamma_s & \delta_s
  \end{pmatrix}
  \begin{pmatrix}
    \xone \\ \yone 
  \end{pmatrix},
\]
for some coefficients $\alpha_s, \beta_s,\gamma_s, \delta_s$ in
$W_R^{n_t}$. 
It follows from the same computation (\ref{eq:46}) of (\ref{46}) that 

\begin{align}
  \label{eq:672}
  \begin{pmatrix}
    \alpha_{s-1} & \beta_{s-1}  \\
    \gamma_{s-1} & \delta_{s-1}
  \end{pmatrix}&=
   \begin{pmatrix}
     U_{s-1} & 1 \\
     U_{s-1} V_{s-1} +\pi^{e\ell_{s-1}} & V_{s-1} 
   \end{pmatrix}
  \begin{pmatrix}
    \alpha_s & \beta_s  \\
    \gamma_s & \delta_s
  \end{pmatrix}, \\   
  \begin{pmatrix}
    \alpha_{t} & \beta_{t}  \\
    \gamma_{t} & \delta_{t}
  \end{pmatrix}&=
   \begin{pmatrix}
     U_t & 1 \\
     U_t V_t +\pi^{e\ell_t} & V_t 
   \end{pmatrix},
\end{align}
where 
\[ U_{s}:=T_{n_s}^{(\ell_{s+1}+\ell_{s+2}+\dots+\ell_t)}\quad
\text{and}\quad V_{s}:=\sum_{i=1}^{\ell_s-1}
T^{(\ell_s+\dots+ \ell_t-i)}_{n_{s-1}+i} \pi ^{ei}. \]
Therefore, 
\begin{equation}
  \label{eq:673}
  \begin{pmatrix}
    \alpha & \beta  \\
    \gamma & \delta
  \end{pmatrix}:=   
  \begin{pmatrix}
    \alpha_1 & \beta_1  \\
    \gamma_1 & \delta_1
  \end{pmatrix}=   
  \prod_{i=1}^{t} 
  \begin{pmatrix}
     U_i & 1 \\
     U_i V_i +\pi^{e\ell_i} & V_i 
  \end{pmatrix}
\end{equation}

\begin{thm}\label{68}
  If $\ul a=(a^i)$ is spaced, then the points in $\scrM_{\ge \ul a}$
  with slope sequence $s(|\ul a|)$ are dense in $\scrM_{\ge \ul a}$.
\end{thm}
\begin{proof}
  We will show that for each point $x=(A_0,\lambda_0, \iota_0)\in
  \scrM_{\ul a}(k)$, there is a deformation in $\scrM_{\ge \ul a}$
  whose generic point has 
  slope sequence $s(|\ul a|)$. Let $M_0$ be the \dieu module of $A_0$. Let
  $R=k[\![t_{i, a^i}]\!], i\in\tau(M_0)$. We construct a deformation
  $M_R$ by (\ref{62}) with $T_{i,j}=0$ except for $i\in \tau(M_0) $
  and $j=a^i$. Note that $V_i=0$ in (\ref{eq:673}). We have (cf. \ref{48})
  \[    \begin{pmatrix}
  F^f \pi^e \xone \\ F^f \yone 
  \end{pmatrix}=\pi^{|\ul a|}A \begin{pmatrix}
   \pi^e \xone \\ \yone 
  \end{pmatrix}, \quad A=\prod_{i=1}^t 
  \begin{pmatrix}
    U_{i} & \pi^{e-a^{n_i}} \\
  \pi^{e(\ell_{i}-1)-a^{n_i}} & 0
  \end{pmatrix}, \]
where $U_i=(c_{n_i}+T_{n_i,a^{n_i}})^{(\ell_{i+1}+\ell_{i+2} 
+\dots+\ell_t)}$ for some $c_{n_i}$ in $W^{n_i}$. We may assume that
$e(\ell_i-1)>a^{n_i}$ for some $i$, otherwise $|\ul a|=\frac{g}{2}$
and there is nothing to prove. Assume that  $e(\ell_1-1)>a^{n_1}$, we
have
\begin{equation}
  \begin{split}
   A&\equiv 
  \begin{pmatrix}
    U_{1} & 0 \\
    0 & 0
  \end{pmatrix}
  \begin{pmatrix}  
    U_{2} & * \\
    * & 0
  \end{pmatrix}\dots
  \begin{pmatrix}
    U_{t} & * \\
    * & 0
  \end{pmatrix} \mod \pi\\ 
  &\equiv 
  \begin{pmatrix}
    P & * \\
    0 & 0
  \end{pmatrix} \mod \pi
  \end{split}
\end{equation}
where $P$ is a nonzero polynomial in $U_i$'s. Hence $P$ is a unit
in $W(K)[\pi]$, where $K$ is the perfection of $\Frac(R)$. Therefore
$M_K=\Cart(K)\otimes_{\Cart(R)} M_R$ has slope sequence $s(|\ul a|)$.\qed
\end{proof}

\begin{cor}\label{69}
  The generic point of each irreducible component of $\scrM_{\ge \ul
  a}$ has slope sequence $\ge s(\lambda(\ul a))$, where $\lambda(\ul a):=\max
  \{|\ul b|; \ul b\le \ul a, \text{$\ul b$ is spaced} \}$.
\end{cor}
\begin{proof}
  This follows from Proposition~\ref{48} and Grothendieck's
  specialization theorem.\qed
\end{proof}

\begin{prop}\label{610}
  Let $\ul a=(a^i)$ be an $a$-type with $a^i\le [e/2]$ for all
  $\iz$. Let $x=(A_0,\lambda_0,\iota_0)\in
  \scrM(k)$ such that $A_0$ is superspecial. Then in every open
  neighborhood of x in $\scrM_{\ge \ul a}$ there exists a point of
  slope sequence $s(|\ul a|)$. 
\end{prop}
\begin{proof}
  Let $c:=[e/2]$ and $M_0$ be the \dieu module of $A_0$. By
  Lemma~\ref{45} (2), we can choose a basis $\{X_i,Y_i\}$ for $M^i$
  such that 
\[ 
\begin{pmatrix}
  FX_{i-1} \\ FY_{i-1} 
\end{pmatrix}=
\begin{pmatrix}
  0 & 1 \\ \pi^e & 0 \\
\end{pmatrix}
\begin{pmatrix}
  X_i \\ Y_i
\end{pmatrix},\quad \forall\, \iz.  \]
Let $R=k[\![t_{i}]\!]_{\iz}$. We construct a deformation
  $M_R$ of $M_0$ by (\ref{62}):
\[  \begin{pmatrix}
  FX_{i-1} \\ FY_{i-1} 
\end{pmatrix}=
\begin{pmatrix}
  \pi^{a^i}T_i & 1 \\ \pi^e & 0 \\
\end{pmatrix}
\begin{pmatrix}
  X_i \\ Y_i
\end{pmatrix},\quad \forall\, \iz,  \]
where $T_i$ is the Teichm{\"u}ller lift of $t_i$. 
We have (cf. \ref{48}, \ref{610})
\[  \begin{pmatrix}
  F^f \pi^c X_{f-1} \\ F^f Y_{f-1} 
  \end{pmatrix}=\pi^{|\ul a|}A \begin{pmatrix}
   \pi^c X_{f-1} \\ Y_{f-1} 
  \end{pmatrix}, \quad A=\prod_{i=0}^{f-1} 
  \begin{pmatrix}
    U_{i} & \pi^{c-a^{i}} \\
  \pi^{e-c-a^{i}} & 0
  \end{pmatrix}, \]
where $U_i=T_{i}^{(f-1-i)}$. We may assume that
$a^i<e/2$ for some $i$, otherwise $|\ul a|=\frac{g}{2}$
and there is nothing to prove. Say $a^0< e/2$, we
have
\begin{equation}
  \begin{split}
   A&\equiv 
  \begin{pmatrix}
    U_{0} & * \\
    0 & 0
  \end{pmatrix}
  \begin{pmatrix}  
    U_{1} & * \\
    * & 0
  \end{pmatrix}\dots
  \begin{pmatrix}
    U_{f-1} & * \\
    * & 0
  \end{pmatrix} \mod \pi\\ 
  &\equiv 
  \begin{pmatrix}
    P & * \\
    0 & 0
  \end{pmatrix} \mod \pi
  \end{split}
\end{equation}
where $P$ is a nonzero polynomial in $U_i$'s. Hence $P$ is a unit
in $W(K)[\pi]$, where $K$ is the perfection of $\Frac(R)$. Therefore
$M_K=\Cart(K)\otimes_{\Cart(R)} M_R$ has slope sequence $s(|\ul a|)$.\qed
\end{proof}
\begin{rem}
  Goren and Oort showed [GO, Thm. 5.4.11] that when $p$ is inert in
  $\FF$, the inequality $\ge$ in Corollary~\ref{69} can be strengthened by
  equality $=$. However, the equality does not hold in general by
  Proposition~\ref{610}. It will be interesting to have the sharp
  formula of the slope sequence of the generic points of {\it any}
  alpha stratum. When the $a$-types are spaced, the equality can be
  achieved as it stands in Theorem~\ref{68}. 
\end{rem}

\subsection{}
\label{611}
Let $x=(A_0, \lambda_0,\iota_0)$ and $M_0$ be as in (\ref{62}). Suppose
that \emph{the reduced $a$-number $t:=|\tau|$ of $x$ is one},
where $\tau$ is the $a$-index of $M_0$. We assume that $\tau=\{0\}$ for
simplicity. By Proposition~4.4, we choose a basis for $M_0$ as in
(\ref{eq:620}) of (\ref{62}).  

Let $\ell=\ord_\pi (c_0)$. We have $a^0=\min\{e,\ell\}$ and
$\slope(M_0)=s(\min\{\frac{g}{2},\ell\})$ by (\ref{49}).  

Take $\ul a=(0,0,\dots 0)$ and construct the universal deformation
$M_R$ of $M_0$ for $\scrM$ by (\ref{62}). It is the \dieu module of the formal
group attached to the universal formal deformation $(\wt A, \wt \lambda, \wt
\iota)$ over $\scrM_x^{\wedge}=\Spf R$. For each $m\in S(g)$, $m\le
\ell$, let $\scrM^{\ge s(m)}_x$ denote the reduced closed
subscheme (in $\scrM_x^{\wedge}$) consisting of points with slope
sequence $\ge s(m)$ (\ref{112}).  
We will find the defining equations of the  
subscheme $\scrM^{\ge s(m)}_x$ and show that the generic point of
$\scrM^{\ge s(m)}_x$ has slope sequence $s(m)$. 

From (\ref{eq:673}), we have
\[  \begin{pmatrix}
  F^f X_0 \\ F^f Y_0 
\end{pmatrix}=
\begin{pmatrix}
  U_1 & 1 \\ U_1 V_1+\pi^g & V_1 
\end{pmatrix}
\begin{pmatrix}
  X_0 \\ Y_0 
\end{pmatrix},  
\] 
where 
\[ U_1=T_0\quad \text{and}\quad V_1=\sum_{i=1}^{f-1}T_i^{(f-i)}\pi^{ei}. \]
Therefore we have the Cayley-Hamilton equation $F^{2f}
X_0-(U_1^{(f)}+V_1)F^f X_0-\pi^g=0$. The subscheme $\scrM^{\ge s(m)}_x$ is
defined by the equations obtained from ``$\ord_\pi(U_1^{(f)}+V_1)\ge m$''.

From (\ref{eq:651}) of (\ref{65}), write 
\begin{equation}
  \begin{split}
    U_1^{(f)}+V_1&=c_0^{(f)}+\sum_{i=0}^{f-1}\sum_{j=0}^{e-1}
    T_{i,j}^{(f-i)}\pi^{ei+j} \\
    &=c_0^{(f)}+\sum_{i=0}^{g-1}T_{k}\pi^{k} 
  \end{split}
\end{equation}
where $T_k:=T_{i,j}^{(f-i)}, t_k=t_{i,j}$ for $k=ei+j, 0\le j\le
e-1$. We can see 
that the defining equations for $\scrM^{\ge s(m)}_x$ 
are $t_0=t_1=\dots t_{m-1}=0$. Let $K_m$ be the perfection of the
generic residue field of $\scrM^{\ge s(m)}_x$. The element $T_m$ is
a unit in $W(K_m)$, therefore the generic point of $\scrM^{\ge s(m)}_x$
has slope sequence $s(m)$.

\begin{thm}\label{612}
  Let $x:\Spec k\to \scrM$ be a geometric point of reduced $a$-number
  one. Then each closed subscheme $\scrM^{\ge s(m)}_x$, for $m\in S(g),\
  s(m)\le \slope(x)$ (\ref{321}), is formally smooth of codimension
  $\lceil m \rceil$ and its generic point has slope sequence $s(m)$, where
  $\lceil m\rceil$ denotes the smallest integer not less than $m$. 
\end{thm}
\begin{cor}\label{613}
  Let $U$ be the subset of $\scrM$ consisting of points with reduced
  $a$-number $\le 1$.  Then 
  the strong Grothendieck conjecture holds for $U$ (\ref{113}).
\end{cor}

%Note that a point in $\scrM\otimes k(v)$ of reduced $a$-number one
%automatically satisfies the Rapoport condition.

\begin{cor}\label{614}
  The strong Grothendieck conjecture for $\scrM$  holds when $p$ is
  totally ramified in $\FF$ (\ref{113}).
\end{cor}
\begin{cor}\label{615}
  The weak Grothendieck conjecture for $\scrM$ holds (\ref{113}).
\end{cor}
\begin{proof}
  It follows from Lemma~{\ref{73}} that there is a supersingular point
  of reduced $a$-number one. Then the assertion follows from
  Theorem~\ref{612}. \qed
\end{proof}

\subsection{}
\label{616}
  \emph{In the rest of this section we assume that $p$ is totally
  ramified in $\FF$}. Denote by $\scrM_v^{\DP}$ the reduction
  $\scrM^{\DP}\otimes_\Z k(v)$ of $\scrM^{\DP}$ modulo $v$. Let
  $x=(A_0,\lambda_0, \iota_0)\in
  \scrM_v^{\DP}(k)$ be a geometric point and let $\lie(A_0)=\{e_1,
  e_2\}$ and $\ul a(A_0)=\{a_1, a_2\}$. We assume that $e_1\le e_2$ and
  $a_1\le a_2$. Let $M_0$ be the covariant \dieu module of $A_0$. We
  can choose two $W[\pi]=\O\otimes W$-bases $\{X_1, X_2\}$,
  $\{Y_1,Y_2\}$ of $M_0$ such that 
  \[ VuY_1=\pi^{e_2}X_1\quad\text{and}\quad  VuY_2=\pi^{e_1}X_2, \] 
  where $pu=\pi^e$. We have 
  \[ FX_1=\pi^{e_1}Y_1\quad \text{and}\quad  FX_2=\pi^{e_2}Y_2. \]
  Write $Y_1=\alpha X_1+\beta X_2$ and $Y_2=\gamma X_1+\delta
  X_2$. Then we have
  \[ 
  \begin{pmatrix}
    FX_1 \\ FX_2
  \end{pmatrix}=
  \begin{pmatrix}
    \pi^{e_1} \alpha & \pi^{e_1}\beta \\ \pi^{e_2} \gamma & \pi^{e_2}
    \delta 
  \end{pmatrix}
  \begin{pmatrix}
    X_1 \\ X_2
  \end{pmatrix}\quad \text{and} \quad
  (F,V)M=<\pi^{e_1}X_2,\pi^{e_2}X_1, \pi^{e_1} \alpha
  X_1>_{W[\pi]}. \]   
  It follows that $a_1=e_1$ and $a_2=\min\{e_2,
  e_1+\ord_\pi(\alpha)\}$. As $e_2\ge \frac{e}{2}$, by the argument in
  (\ref{410}) we have $\slope(A_0)=s(i)$, where $i=\min\{\frac{e}{2},
  e_1+\ord_\pi(\alpha)\}$. Note that the $a$-type $\ul a(A_0)$
  determines the invariant $\lie(A_0)$ and $\slope(A_0)$: $e_1=a_1$
  and $\slope(A_0)=s(\min\{\frac{e}{2}, a_2\})$.

\subsection{}
\label{617}
In [DP], Deligne and Pappas have defined a closed algebraic substack
$\scrN_i$ of $\scrM_v^{\DP}$ which classifies the objects of Lie type
$\{e_1, e_2\}$ with $e_1\ge i$. They have shown that the complement
$\scrN^o_i$ of $\scrN_{i+1}$ in $\scrN_i$ is a smooth algebraic stack
of dimension $e-2i$ if it is non-empty. It follows from their results
and Theorem~\ref{75} that points of Lie type $\{e_1,e_2\}$ are dense
in $\scrN_{e_1}$. It follows from (\ref{616}) that any point in
$\scrN_{e_1}$ has slope sequence $\ge s(e_1)$. The following lemma
confirms the density of points with slope sequence $s(e_1)$ 
in $\scrN_{e_1}$. 

\begin{lemma}\label{618}
  If $a_2>a_1$, then there is a deformation $(\wt A, \wt \lambda, \wt
  \iota)$ over $k[\![t]\!]$ of $(A_0,\lambda_0,\iota_0)$ whose generic
  point has $a$-type $(a_1, a_2-1)$.
\end{lemma}
\begin{proof}
  As the forgetful map $\Def[A_0,\lambda_0,\iota_0]\to
  \Def[A_0,\iota_0]$ induces an equivalence of deformation functors in
  $\scrN_{e_1}$, we will construct a deformation of abelian
 $O_\FF$-varieties in $\scrN_{e_1}$. By the reduction step in [DP, 4.3] and
  the construction of (\ref{62}), we can construct a \dieu $\O$-module
  $M_R$ over $R:=k[\![t]\!]$ of $M_0$ such that
  \begin{equation}
    \label{eq:6181}
    \begin{split}
      FX_1&=\pi^{e_1}\alpha X_1+\pi^{e_1}\beta(X_2+T \pi^{a_2-e_1-1}
      X_1)\\
      FX_2&=\pi^{e_2}\gamma X_1+\pi^{e_2}\delta(X_2+T \pi^{a_2-e_1-1}
      X_1),\\ 
    \end{split}
  \end{equation}
  where $T$ is the Teichm{\"u}ller lift of $t$. Note that $\beta$ is a
  unit in $W[\pi]$: if $\beta\equiv 0 \mod \pi$, then $\alpha$ is a
  unit and $a_2=a_1$.
  By base change to $K:=k((t))^{\mathrm{perf}}$, we have 
  \[ (F,V)M_K=<\pi^{e_1}X_2,\pi^{e_2}X_1, 
  (\pi^{e_1} \alpha+\pi^{a_2-1}\beta T) X_1>_{W(K)[\pi]}, \]
  thus $\ul a(M_K)=(a_1, a_2-1)$. This completes the proof. \qed 
\end{proof}

\begin{thm}\label{619}
  The strong Grothendieck conjecture for $\scrM_v^{\DP}$ holds when
  $p$ is totally ramified in $\FF$.
\end{thm}
\begin{proof}
  Let $x\in \scrM_v^{\DP}(k)$ be a geometric point with $a$-type
  $(a_1,a_2)$ and slope sequence $s(i)$, where $i=\min\{\frac{e}{2},a_2\}$. It
  suffices to deform the point $x$ to a point with slope sequence
  $s(i-1)$. If $a_2>a_1$, then by Lemma~\ref{618} we can deform to a
  point with slope sequence $s(i-1)$. Suppose that $a_2=a_1=e_1$, we have
  $\slope(x)=s(e_1)$. By (\ref{617}), we can deform to a point of Lie
  type $\{e_1-1, e_2+1\}$. As points with slope sequence $s(e_1-1)$ 
  are dense in $\scrN_{e_1-1}$, we can further deform to a point with
  slope sequence $s(e_1-1)$. \qed 
\end{proof}

\begin{cor}\label{620}
  Each Newton stratum of $\scrM_v^{\DP}$ with slope sequence $s(m)$,
  $m\in S(g)$ (\ref{321}), has pure
  dimension $g-\lceil m\rceil$ when $p$ is totally ramified in
  $\FF$. 
\end{cor}
\begin{proof}
  This follows from the purity of Newton strata [dJO] and Theorem~\ref{619}.
\end{proof}

\begin{rem} 
  In [C] Chai gives a group-theoretic dimension formula for Newton
  strata arising from quasi-split groups. He expects that it is so for
  good reduction of PEL-type Shimura varieties [C, Question 
  7.6, p.~984]. Corollary~\ref{620} suggests that his description can
  be applied for a larger class, not necessarily restricted to the
  good reduction case.  
\end{rem}

\section{An algebraization theorem}
\label{sec:07}

\subsection{}
Let notations be as in (\ref{11}) and (\ref{112}). 
The set of possible Newton polygons in question is parameterized by
$S(g_1)\times \dots\times S(g_s)$ (\ref{111}), 
where $g_i=[\FF_{v_i}:\Q_p]$. For each
abelian $O_\FF$-variety $A$, the associated $p$-divisible group
$A(p):=A[p^\infty]$ has a decomposition $A(p)=A(p)_{v_1}\oplus \dots
\oplus A(p)_{v_s}$.

Recall (\ref{15}) that a quasi-polarized $p$-divisible $\O$-group
$(H,\lambda,\iota)$ over $k$ is \emph{algebraizable} if it is 
attached to a polarized abelian $O_\FF$-variety 
$(A, \lambda_A, \iota_A)$ over $k$.  

\begin{lemma}\label{73}
  Any supersingular quasi-polarized $p$-divisible $\O$-group $(H,\lambda,
  \iota)$ over $k$ is algebraizable.  
\end{lemma}
\begin{proof}
  Let $E$ be a supersingular elliptic curve over $k$, and
  $A':=E\otimes O_\FF$. It 
  is clear that $A'$ satisfies the Rapoport condition. By [R, Prop.~1.10],
  there exists a separable $O_\FF$-linear polarization $\lambda'$ on
  $A'$. Let $(H_1, \lambda_1, \iota_1)$ be the $p$-divisible group attached
  to $(A',\lambda',\iota')$. It follows from Corollary ~\ref{36} that $(H_1,
  \lambda_1, \iota_1)$ is isogenous to $(H, \lambda,\iota)$. By a
  theorem of Tate, there is a polarized abelian $O_\FF$-variety $(A,
  \lambda_A, \iota_A)$ whose $p$-divisible group is isomorphic to
  $(H,\lambda,\iota)$. \qed 
\end{proof}

\begin{thm}\label{74}
(1) The weak Grothendieck conjecture for $\scrM$ holds. 

(2) The strong Grothendieck conjecture for $\scrM^{\DP}_p$ holds when
    all residue degrees $f_i$ are one.
\end{thm}
\begin{proof}
  (1) It follows from Lemma~\ref{73} that there is a supersingular
      polarized abelian $O_\FF$-variety $A$ such that each component
      $A(p)_{v_i}$ of the associated $p$-divisible group $A(p)$ has
      reduced $a$-number one. Then the assertion follows from the
      theorem of Serre-Tate and Theorem~\ref{612}.

  (2) This follows from the theorem of Serre-Tate and
      Theorem~\ref{619}. \qed
\end{proof}

\begin{thm}\label{75}
  Any quasi-polarized $p$-divisible $\O$-group $(H,\lambda,\iota)$
  over $k$ is algebraizable.
\end{thm}
\begin{proof}
  It follows from  Theorem~\ref{74} (1) that
  any slope sequence in $S(g_1)\times \dots\times S(g_s)$ can be realized
  by a point in $\scrM$. Then the theorem follows from
  Corollary~\ref{36}. \qed
\end{proof}

\section{An example}
\label{sec:08}

In this section we construct a separably polarized abelian $O_\FF$-scheme
$A$ over a complete DVR $R$ whose close fibre $A_k$ does not satisfy
the Rapoport condition. 

Let $g=2$ and $p> 3$ be a ramified prime in the totally quadratic
real field $\FF$. We have $O_\FF\otimes W(k)=W(k)[\pi], \pi^2=p$. Let $M$
be a free $W(k)[\pi]$-module generated by $e_1, e_2$ with the Verschiebung
action
\[ Ve_1=\pi e_2, Ve_2=\pi e_1. \]
and with the alternating form determined by
\[ \<e_1, e_2'\>=\< e_1',e_2\>=1 \]
and other pairing are $0$ for the $W$-basis $e_1, e_1', e_2,
e_2'$, where $e_1':=\pi e_1, e_2':=\pi e_2$. 
By the algebraization theorem, there is a polarized abelian
$O_\FF$-variety $A_0$ over $k$ with the prescribed \dieu module $M$.

We have the Hodge filtration of $M/pM=H^1_{DP}(A_0/k)^*$
\[ 0\to VM/pM=k\bar e_1'\oplus k \bar e_2'\to M/pM \to M/VM=k\bar e_1
\oplus\bar e_2\to 0. \]
Clearly, $A_0$ does not satisfy the Rapoport condition, as 
$\pi=0$ on $M/VM$.  

Let $R:=W(k)[\sqrt{p}]$. By Grothendieck-Messing's
theory, Serre-Tate's Theorem and Grothendieck's Existence Theorem, we
can lift the abelian variety $A_0$ with the additional structure over
$R$ by lifting the 
Hodge filtration with respect to the addition structure, see [Y2,
Sect. 4]. 
Let $N$ be the $R$-submodule of $M\otimes_W R$ generated by
$e_1'+\sqrt{p}e_1$ and $e_2'-\sqrt{p}e_2$. It is easy to check that
$N$ is stable by $O_\FF$-action, $N\otimes_R k=VM/pM$ and
$\<N,N\>=0$. Thus, we get a desired polarized abelian $O_\FF$-scheme
over $R$.

\section{A computation of the Hecke correspondence}
\label{sec:09}

\subsection{}

Let $M_0$ be the quasi-polarized \dieu $\O$-module of the polarized
abelian variety $A_0$ in the previous
section. Let $N$ be the quasi-polarized \dieu module containing $M_0$
with $VN=M$ and $\<\, , \>_N=\frac{1}{p} \<\, , \>$. We have 
\[ \begin{array}{ll}
    Ve_1=\pi e_2=e_2'  & \quad Fe_2'=pe_1 \\
    Ve_2=\pi e_1=e_1'  & \quad Fe_1'=pe_2 \\
    Ve_1'=\pi Ve_1=pe_ 2 & \quad Fe_2=e_1' \\
    Ve_2'=\pi Ve_2=pe_1  & \quad Fe_1=e_2' \\
\end{array} \]
and $\<e_1,e_2'\>=1, \<e_1',e_2\>=1$. From $N=\frac{1}{p}FM$, we have 
\[ N=<e_1, e_2, \frac{1}{p}e_1',  \frac{1}{p}e_2'>_W. \]
Write $X_1=e_1, X_2=e_2, X_1'=\frac{1}{p}e_1', X_2'=\frac{1}{p}e_2'$, we have
\[ M_0=<X_1,X_2,pX_1',pX_2'>_W. \]
and $\<X_1,X_2'\>_N=\<X_1',X_2'\>_N=1$. Denote by $\bar N$ the
quotient of $N$ modulo $pN$ and $\bar M_0$ the image of $M_0$ in $\bar
N$. Write $x_i,x'_i$ the image of $X_i, X'_i$ in $\bar N$.

\subsection{}
Let $\calX$ be space of $\pi$-invariant maximal isotropic ``\dieu'' subspaces
of $\bar N$ over $k$
\[ \calX(k)=\{ 0\subset \bar M \subset \bar N; \dim \bar M=2, \ 
\pi \bar M\subset \bar M, F\bar M\subset \bar M, V\bar M\subset \bar
M, \<\bar M,\bar M\>_N =0 \} \]
We regard $\calX$ as a reduced subscheme over $k$. It is a closed
subscheme of the Lagrangian Grassmanian $\mathrm LG(2,4)$, hence a
projective variety. 

There is a finite morphism $\mathrm{pr}:\calX\to \scrM^{\DP}_p$
(cf. [Y1. Sect. 6])
sending $\bar M_0 \mapsto (A_0,\lambda_0,\iota_o)$ and the morphism
$\mathrm{pr}$ factors through the supersingular locus $\scrS^{\DP}$ of
$\scrM^{\DP}_p$.

We have $\bar{M}_0=<x_1, x_2>_k\in \calX(k)$. Let
$\bar M_t\in \calX$ be the $k$-subspace of $\bar N$ generated by
\[ \tilde x_1=x_1+t_{11}x'_1+t_{12}x'_2, \quad \tilde
x_2=x_2+t_{21}x'_1+t_{22}x'_2. \] 
The points $\bar M_t$ form a Zariski open neighborhood of $\bar M_0$
in $\calX$, which we denote by $\calU$. We will show that 
$\calU\simeq
\Spec k[t_1,t_2,t_3]/(t_1^{p+1}-t_2^{p+1}, t_1^2+t_2t_3)$. 

\subsection{}

From $\<\bar M_t, \bar M_t\>_N=0$, we get $t_{11}+t_{22}=0$. 
From $X'_i=\frac{1}{p}e'_i=\frac{1}{\pi}e_i=\frac{1}{\pi}X_i$, we have
\[ \pi x'_i=x_i,\  \pi x_i=0, i=1,2. \]
One computes
\begin{equation*}
  \begin{split}
    \pi \tilde x_1&=t_{11}x_1+t_{12}x_2
                  =t_{11}\tilde x_1+t_{12}\tilde
                  x_2-(t_{11}^2+t_{12}t_{21})x'_1 
                  -(t_{11}t_{12}+t_{12}t_{22})x'_2  \\
  \end{split}
\end{equation*}
and concludes from $\pi \tilde x_1\in \bar M_t$ that 
\begin{equation}
  \label{eq:831}
  t_{11}^2+t_{12}t_{21}=0, \ \text{and } t_{12}(t_{11}+t_{22})=0.
\end{equation}
Similarly from $\pi \tilde x_2\in \bar M_t$, one gets
\begin{equation}
  \label{eq:832}
  t_{22}^2+t_{12}t_{21}=0, \ \text{and } t_{21}(t_{11}+t_{22})=0.
\end{equation}

For the stability of $\bar M_t$ by $F$ and $V$, one computes
\[ F \tilde x_1=t_{11}^px_2+t_{12}^px_1
                  =t_{11}^p\tilde x_2+t_{12}^p\tilde
                  x_1-(t_{11}^p t_{21}+t_{12}^p t_{11})x'_1 
                  -(t_{11}^pt_{22}+t_{12}^{p+1})x'_2. \]
and obtains 
\begin{equation}
  \label{eq:833}
  t_{11}^p t_{21}+t_{12}^p t_{11}=0
\end{equation}
\begin{equation}
  \label{eq:834}
  t_{11}^pt_{22}+t_{12}^{p+1}=0
\end{equation}
Similarly from $F \tilde x_2\in \bar M_t$, one obtains
\begin{equation}
  \label{eq:835}
  t_{21}^{p+1}+t_{22}^p t_{11}=0
\end{equation}
\begin{equation}
  \label{eq:836}
  t_{21}^pt_{22}+t_{22}^p t_{12}=0
\end{equation}

Applying $V\bar M_t \subset \bar M_t$, one does not get new
equations but (\ref{eq:833})--(\ref{eq:836}). From (\ref{eq:834}), one has
$t_{12}=\alpha t_{11}, \alpha^{p+1}=1$. From (\ref{eq:835}), one has
$t_{21}=\beta t_{22}=-\beta t_{11}, \beta^{p+1}=1$. From (\ref{eq:831}),
one has $1=\alpha\beta$. These parameters satisfy the the equations
(\ref{eq:833}) and (\ref{eq:836}).  

We computed that $\calU=\{(t_{11}, t_{12}, t_{21}, t_{22})=(t,\alpha
t, -\frac{1}{\alpha} t, -t); \alpha^{p+1}=1 \}$, hence \\
$\calU\simeq\Spec k[t_1,t_2,t_3]/(t_1^{p+1}-t_2^{p+1},
t_1^2+t_2t_3)$. 
Compared with a result of [BG, p.~476, 3.], the morphism $\mathrm{pr}$ maps
$\calX$ \emph{onto} the ($p+1$) irreducible components of
$\scrS^{\DP}$ containing $(A_0,\lambda_0,\iota_0)$. 

\begin{prop}
  Let $\calU$ be as above. 
  There is a Zariski open neighborhood $\scrV$ of $0$ in $\calU$ such
  that $\mathrm{pr}: \scrV\to \scrS^{\DP}$ is an
  {\'e}tale neighborhood of $\scrS^{\DP}$ at $(A_0,\lambda_0,\iota_0)$.
\end{prop}
\begin{proof}
  Choose a finite {\'e}tale cover $\scrM^{\DP}_{p,n}\to \scrM^{\DP}_{p}$
  by adding a prime-to-$p$ level structure with $n\ge 3$. Choose a
  lift $(A_0,\lambda_0,\iota_0,\eta_0)$ of
  $(A_0,\lambda_0,\iota_0)$. Then there is a lift 
  $\wt{\mathrm{pr}}:\calX\to\scrM^{\DP}_{p,n}$ which sends $\bar{M_0}$ to
  $(A_0,\lambda_0,\iota_0,\eta_0)$. The morphism $\wt{\mathrm{pr}}$ becomes a
  closed immersion as the automorphisms of the objects are trivial. It again
  factors through the supersingular locus $\scrS^{\DP}_n$. As $\calX$
  and $\scrS^{\DP}_n$ are reduced schemes, $\calX$ is isomorphic to its
  image. The image is the union of the $p+1$ irreducible components
  of $\scrS^{\DP}$ containing
  $(A_0,\lambda_0,\iota_0,\eta_0)$. Then there is a Zariski open
  $\scrV$ of $\calU$ such that the $\scrV$ is an
  {\'e}tale neighborhood of $\scrS^{\DP}$ at $(A_0,\lambda_0,\iota_0)$. \qed
\end{proof}

\end{document}